\input amstex
\documentstyle {amsppt}
\magnification=1200  
\vsize=9.5truein 
\hsize=6.5truein 
\nopagenumbers 
\nologo

\def\abs#1{\left\vert #1 \right\vert}

\topmatter

\centerline{26 ao\^ut 2005}

\title
Groupes fondamentaux des vari\'et\'es de dimension $3$
\\
et alg\`ebres d'op\'erateurs
\endtitle

\rightheadtext{Groupes de $3$-vari\'et\'es et alg\`ebres
d'op\'erateurs}

\author
Pierre de la Harpe et Jean-Philippe Pr\'eaux
\endauthor


\subjclass
{ Primary 57M05, 57M27. 
Secondary 20F34, 46L35}
\endsubjclass

\keywords 
Manifold of dimension $3$, Seifert manifold,
group with infinite conjugacy classes, factor of type $II_1$.
Vari\'et\'e de dimension $3$, vari\'et\'e de Seifert,
groupe \`a classes de conjugaison infinies, facteur de type $II_1$
\endkeywords

\thanks
Les auteurs remercient le {\it Fonds national suisse de la
recherche scientifique} pour son soutien.
Le second auteur a \'et\'e partiellement financ\'e par le FNSRS.
\endthanks

\abstract\nofrills{ABSTRACT.}
We provide a geometric characterization of manifolds of dimension $3$
with fundamental groups of which all conjugacy classes except $1$ are infinite,
namely of which the von~Neumann algebras are factors of type $II_1$:
they are essentially the $3$-manifolds with infinite fundamental groups
on which there does not exist any Seifert fibration.

Otherwise said and more precisely,
let $M$ be a compact connected $3$-manifold 
and let $\Gamma$ be its fundamental group,
supposed to be infinite and with at least 
one finite conjugacy class besides $1$.
If $M$ is orientable, then $\Gamma$ is the fundamental group 
of a Seifert manifold;
if $M$ is not orientable, then $\Gamma$ is the fundamental group 
of a Seifert manifold modulo $\Bbb P$ 
in the sense of Heil and Whitten \cite{HeWh--94}.

We make heavy use of results on $3$-manifolds,
as well classical results 
(as can be found in the books of Hempel, Jaco, and Shalen),
as more recent ones (solution of the Seifert fibred space conjecture).

\medskip

\noindent
{\rm R\'ESUM\'E.} 
Nous proposons une caract\'erisation g\'eom\'etrique
des vari\'et\'es de dimension~$3$
ayant des groupes fondamentaux
dont toutes les classes de conjugaison autres que~$1$ sont infinies,
c'est-\`a-dire dont les alg\`ebres de von Neumann 
sont des facteurs de type~$II_1$~:
ce sont essentiellement
les $3$-vari\'et\'es \`a groupes fondamentaux infinis
qui n'admettent pas de fibration de Seifert.

Autrement dit et plus pr\'ecis\'ement, 
soient~$M$ une $3$-vari\'et\'e connexe compacte
et $\Gamma$ son groupe fondamental,
qu'on suppose \^etre infini 
et avec au moins une classe de conjugaison finie autre que~$1$.
Si $M$ est orientable, alors~$\Gamma$ est groupe fondamental 
d'une vari\'et\'e de Seifert~;
si $M$ est non orientable, alors $\Gamma$ est groupe fondamental 
d'une vari\'et\'e de Seifert modulo $\Bbb P$ 
au sens de Heil et Whitten \cite{HeWh--94}.

Nous faisons un usage intensif de r\'esultats concernant les $3$-vari\'et\'es,
autant classiques 
(comme on les trouve dans les livres de Hempel, Jaco et Shalen)  
que plus r\'ecents
(solution de la  conjecture des fibr\'es  de Seifert).

\medskip

\noindent
I.~Introduction.
\par\noindent
II.~Rappels sur les $3$-vari\'et\'es,
la d\'ecomposition de Kneser-Milnor, 

et le th\'eor\`eme de Grushko-Stallings.
\par\noindent
III.~Les vari\'et\'es de Seifert.
\par\noindent
IV.~Rappels sur
la conjecture des fibr\'es de Seifert (cas orientable)

et sur un r\'esultat de Hempel-Jaco.
\par\noindent
V.~Groupes cci.
\par\noindent
VI.~Le th\'eor\`eme principal pour les vari\'et\'es orientables irr\'eductibles.
\par\noindent
VII.~Le th\'eor\`eme principal pour les vari\'et\'es orientables.
\par\noindent
VIII.~Rappel sur la conjecture des fibr\'es de Seifert (cas non orientable).
\par\noindent
IX.~Le th\'eor\`eme principal 
pour les vari\'et\'es non orientables $\Bbb P^2$-irr\'eductibles.
\par\noindent
X.~Le th\'eor\`eme principal 
pour les vari\'et\'es non orientables.
\par\noindent
XI.~Cas d'un entrelacs et cas d'une vari\'et\'e hyperbolique.
\par\noindent
XII.~Groupes de dimension cohomologique trois \`a dualit\'e de Poincar\'e.
\par\noindent
XIII.~La propri\'et\'e cci forte
--
Une application \`a certaines repr\'esentations unitaires.

\endabstract

\address
Pierre de la Harpe, 
\newline
Section de Math\'ematiques, Universit\'e de
Gen\`eve, C.P. 64,  CH-1211 Gen\`eve 4.
\newline
Mel~: Pierre.delaHarpe\@math.unige.ch
\endaddress

\address
Jean-Philippe Pr\'eaux, 
\newline
Centre de recherche de l'Ecole de l'air,
Ecole de l'Air,  
F-13661 Salon de Provence air. 
\newline
Centre de math\'ematiques et d'informatique,
Universit\'e de Provence,
39 rue F.~Joliot-Curie, 13453 F-Marseille cedex 13.
\newline
Mel~: preaux\@cmi.univ-mrs.fr
\endaddress

\endtopmatter

\document

\head 
{\bf 
I.~Introduction
}
\endhead
\medskip

Plusieurs classes de groupes ont \'et\'e \'etudi\'ees 
du point de vue des alg\`ebres d'op\'erateurs. 
Parmi les plus anciens exemples de la litt\'erature 
figurent les groupes libres,
et plus g\'en\'era\-lement les produits libres~; 
voir par exemple le \S \ 6.2 de \cite{ROIV} 
pour le facteur d\'efini par un groupe libre de rang deux. 
Nous abordons ici de ce point de vue l'\'etude
des {\it $3$-groupes}, 
c'est-\`a-dire (ici !) de groupes fondamentaux 
des vari\'et\'es compactes connexes de dimension trois.

Soit $\Gamma$ un groupe. Soient $\xi,\eta$ deux fonctions \`a
valeurs complexes sur $\Gamma$, l'une au moins \`a support fini~;
rappelons que leur {\it produit de convolution} est d\'efini par
$$
(\xi \ast \eta)(\gamma) = \sum_{\gamma_1,\gamma_2 \in \Gamma,
\gamma_1\gamma_2 = \gamma} \xi(\gamma_1)\eta(\gamma_2) .
$$
L'alg\`ebre $\Bbb C[\Gamma]$ du groupe $\Gamma$ est l'alg\`ebre de
convolution des fonctions \`a supports finis~; 
elle est munie d'une involution d\'efinie par 
$\varphi^* (\gamma) = \overline{\varphi(\gamma^{-1})}$ 
pour tous $\varphi \in \Bbb C[\Gamma]$ et $\gamma \in \Gamma$. 
Notons $\Cal B (\ell^2(\Gamma))$ l'alg\`ebre involutive
des op\'erateurs lin\'eaires born\'es 
sur l'espace de Hilbert $\ell^2(\Gamma)$. 
La repr\'esentation r\'eguli\`ere gauche 
$\lambda_{\Gamma} : \Bbb C[\Gamma] \longrightarrow \Cal B (\ell^2(\Gamma))$ 
est d\'efinie par 
$\left(\lambda_{\Gamma}(\varphi)\right)(\xi) = \varphi \ast \xi$ 
pour tous $\varphi \in \Bbb C [\Gamma]$ et $\xi \in \ell^2(\Gamma)$~; 
c'est une $*$-repr\'esentation fid\`ele. 
{\it L'alg\`ebre de von Neumann} $W^*_{\lambda}(\Gamma)$ de $\Gamma$
est l'adh\'erence pour la topologie faible de
$\lambda_{\Gamma}(\Bbb C[\Gamma])$ dans $\Cal B (\ell^2(\Gamma))$.
La {\it C$^*$-alg\`ebre r\'eduite} $C^*_{\lambda}(\Gamma)$ de
$\Gamma$ est l'adh\'erence pour la topologie normique de
$\lambda_{\Gamma}(\Bbb C[\Gamma])$ dans $\Cal B (\ell^2(\Gamma))$.

Notre but est de d\'egager certaines propri\'et\'es 
de l'alg\`ebre $W^*_{\lambda}(\Gamma)$  
lorsque $\Gamma$ est un $3$-groupe.
(Plus tard, on peut esp\'erer comprendre aussi $C^*_{\lambda}(\Gamma)$).
Les preuves consistent essentiellement \`a recueillir et combiner
des r\'esultats connus concernant les $3$-vari\'et\'es
(voir notamment \cite{Hemp--76}, \cite{Jaco--77} et \cite{JaSh--79}).
En particulier, nous utilisons deux r\'esultats profonds
caract\'erisant certaines $3$-vari\'et\'es dont le groupe fondamental
poss\`ede un sous-groupe distingu\'e de type fini non r\'eduit \`a~$1$~:
la {\it conjecture des fibr\'es de Seifert} 
(aujoud'hui un th\'eor\`eme, voir plus bas), 
qui traite du cas o\`u ce sous-groupe est cyclique, 
et un th\'eor\`eme de Hempel et Jaco,
qui traite de cas o\`u ce sous-groupe est d'indice infini et non cyclique.

Un groupe $\Gamma$ est dit {\it \`a classes de conjugaison infinies}, 
ou plus bri\`evement {\it cci}, s'il est infini et
si ses classes de conjugaison distinctes de $1$ sont toutes infinies.
La caract\'erisation suivante des groupes cci est classique~: 
c'est le lemme 5.3.4 de \cite{ROIV}~; 
voir aussi, par exemple, le lemme 4.2.18 de \cite{Saka--71}. 
Rappelons au pr\'ealable 
qu'une alg\`ebre de von Neumann $M$ est un {\it facteur de type $II_1$} 
si elle poss\`ede les trois propri\'et\'es suivantes~:
\roster
\item"(i)"
le centre de $M$ est r\'eduit \`a $\Bbb C$,
\item"(ii)"
il existe une forme lin\'eaire non nulle $\tau$ sur $M$ qui est une
trace, c'est-\`a-dire telle que $\tau(xy-yx)=0$ pour tous $x,y \in M$,
\item"(iii)"
l'alg\`ebre $M$ est de dimension infinie.
\endroster

\medskip

\proclaim{Caract\'erisation de Murray et von Neumann}
L'alg\`ebre de von Neumann $W^*_{\lambda}(\Gamma)$ 
est un facteur de type $II_1$ si et seulement si 
le groupe $\Gamma$ est cci.
\endproclaim

\medskip

Commen\c cons par un \'enonc\'e 
concernant les groupes fondamentaux de surfaces,
qui est une cons\'equence presqu'imm\'ediate
de la classification des surfaces
(voir aussi la proposition 3 ci-dessous).

\medskip

\proclaim{Cas des $2$-groupes}
Soient $F$ une surface compacte connexe et $\Gamma$ son groupe fondamental
($F$ peut \^etre orientable ou non, \`a bord vide ou non vide).
Le groupe $\Gamma$ est infini et non cci si et seulement si 
$F$ est hom\'eomorphe \`a l'une des surfaces de la liste
\par
\centerline{anneau, ruban de M\"obius, $2$-tore, bouteille de Klein,}
\par\noindent
c'est-\`a-dire si et seulement si $F$ admet un feuilletage par des cercles.
\endproclaim

\medskip

Le but principal du pr\'esent travail 
est de montrer une assertion analogue pour les $3$-groupes.
Le lecteur expert en $3$-vari\'et\'es 
pourra aborder l'essentiel de nos arguments 
en parcourant les preuves des th\'eor\`emes 12 et 13,
consacr\'es au cas orientable.
Pour le cas non orientable, voir le th\'eor\`eme 18 
(qui compl\`ete l'\'enonc\'e ci-dessous).

\medskip

\proclaim{Th\'eor\`eme}
Soient $M$ une $3$-vari\'et\'e et $\Gamma$ son groupe fondamental~;
on suppose que $\Gamma$ est infini et non cci.

Si $M$ est orientable, $\Gamma$ est groupe fondamental d'une vari\'et\'e de
Seifert.
Si $M$ est non orientable, $\Gamma$ poss\`ede un sous-groupe d'indice $2$
qui est groupe fondamental d'une vari\'et\'e de Seifert.
\endproclaim

\bigskip

Nous remercions Michel Boileau, David Epstein  et Claude Weber 
pour plusieurs commentaires utiles \`a ce texte. 
De plus, le second auteur tient \`a remercier les math\'ematiciens genevois, 
et nomm\'ement Goulnara Arzhantseva, pour leur accueil en 2004 et 2005.

\bigskip
\head{\bf 
II.~Rappels sur les $3$-vari\'et\'es,
la d\'ecomposition de Kneser-Milnor, 
et le th\'eor\`eme de Grushko-Stallings
}
\endhead
\medskip

Soit $M$ une {\it $3$-vari\'et\'e}, c'est-\`a-dire ici une
vari\'et\'e de dimension $3$, compacte, connexe,
qui peut \^etre ou bien close 
(c'est-\`a-dire \`a bord $\partial M$ vide)  
ou bien bord\'ee (c'est-\`a-dire \`a bord non vide). 
Nous convenons de penser \`a $M$ et ses sous-vari\'et\'es 
du point de vue p.l., ou \lq\lq presque lin\'eaire\rq\rq , 
et tous les plongements consid\'er\'es sont suppos\'es p.l.~;
mais les points de vue diff\'erentiable et topologique localement plat 
sont \'equivalents (voir le d\'ebut du chapitre~1 de \cite{Hemp--76}).
Nous notons syst\'ema\-tiquement
$$
\Gamma \, = \, \pi_1(M)
$$
le groupe fondamental de $M$~; il est de pr\'esentation finie.

Le cas des $3$-vari\'et\'es \`a groupes fondamentaux finis 
est le sujet de nombreux travaux remontant au moins aux ann\'ees 1920~; 
voir par exemple \cite{Hopf--25}, \cite{Miln--57}
et \cite{Rubi--95}. 
Nous nous int\'eressons ici en premier lieu aux $3$-groupes {\it infinis.}
Pour ce qui suit, voir \cite{Hemp--76}, 
en particulier le chapitre 3 pour le th\'eor\`eme de Kneser-Milnor
et le chapitre 7 pour celui de Gruschko-Stallings.

Etant donn\'e deux $3$-vari\'et\'es $M_1,M_2$, 
deux $3$-boules plong\'ees $B_1 \subset M_1$, $B_2 \subset M_2$
et un hom\'eomorphisme $\varphi$ du bord de $B_1$ sur celui de $B_2$,
on d\'efinit la {\it somme connexe} $M_1 \sharp\, M_2$
(voir \cite{Hemp--76}, chapitre 3). 
Si $M_1$ et $M_2$ sont orient\'ees
et si $\varphi$ est tel que les orientations naturelles
de $\varphi(\partial B_1)$ et $\partial B_2$ sont oppos\'ees,
le type d'hom\'eomorphisme de la somme connexe 
est ind\'ependant des autres choix~;
dans le cas g\'en\'eral, 
il y a au plus deux types d'hom\'eomorphisme pour $M_1 \sharp\, M_2$.
Dans tous les cas,  le groupe fondamental $\pi_1(M_1 \sharp\, M_2)$
est isomorphe au produit libre $\pi_1(M_1) \ast \pi_1(M_2)$.

Une $3$-vari\'et\'e $M$ est {\it ind\'ecomposable} si, 
pour toute somme connexe $M_1 \sharp\, M_2$ hom\'eo\-morphe \`a $M$, 
l'un au moins des facteurs $M_1,M_2$ est une $3$-sph\`ere.
Une sph\`ere plong\'ee dans une $3$-vari\'et\'e est {\it inessentielle}
si elle borde une $3$-boule et {\it essentielle} sinon.
Une $3$-vari\'et\'e $M$ est {\it irr\'eductible} 
si elle ne poss\`ede aucune sph\`ere essentielle.
Une $3$-vari\'et\'e ind\'ecomposable est ou  bien irr\'eductible, 
ou bien un $\Bbb S^2$-fibr\'e sur un cercle
\footnote{
Nous notons $\Bbb S^n$ [respectivement $\Bbb D^n$, $\Bbb P^n$]
la sph\`ere [resp. le disque, l'espace projectif] de dimension $n$,
et $I = [0,1]$ l'intervalle unit\'e de la droite r\'eelle.
Au lieu de \lq\lq $n$-disque\rq\rq , on \'ecrit aussi
\lq\lq $n$-boule\rq\rq \ ou \lq\lq $n$-cellule\rq\rq .
}
(lemme 3.13 de \cite{Hemp--76}). 
L'espace total d'un $\Bbb S^2$-fibr\'e sur un cercle est 
ou bien orientable,
et alors hom\'eomorphe \`a $\Bbb S^2 \times \Bbb S^1$, 
ou bien non orientable, 
et alors hom\'eomorphe \`a l'espace total du
fibr\'e en sph\`eres non trivial de base un cercle,
not\'e $\Bbb S^2 \tilde \times \Bbb S^1$. 

\medskip

\proclaim{D\'ecomposition de Kneser-Milnor} 
{\it Toute $3$-vari\'et\'e s'\'ecrit
comme somme connexe de vari\'et\'es ind\'ecomposables
$$
M \, = \, X_1 \sharp\, \cdots \sharp\, X_k \sharp\, 
          Y_1 \sharp\, \cdots \sharp\, Y_l
$$
o\`u $k,l \ge 0$ sont des entiers, o\`u $X_1,\hdots,X_k$ sont
irr\'eductibles, et o\`u $Y_1,\hdots,Y_l$ sont des fibr\'es en
$2$-sph\`eres sur des cercles 
($M$ est hom\'eomorphe \`a $\Bbb S^3$ si $k+l=0$).

Lorsque $M$ est orientable, cette d\'ecomposition est unique 
\`a l'ordre pr\`es et \`a hom\'eo\-morphisme pr\`es des facteurs.
Lorsque $M$ est non orientable, il existe une unique 
\footnote{
Rappelons toutefois que, pour toute $3$-vari\'et\'e non orientable $M$,
les sommes connexes $M \sharp\, \left(\Bbb S^2 \tilde \times \Bbb S^1\right)$
et $M \sharp\, \left(\Bbb S^2 \times \Bbb S^1\right)$
sont hom\'eomorphes (lemme 3.17 de \cite{Hemp--76}). 
}
telle d\'ecomposition dans laquelle aucun des $Y_j$ 
n'est hom\'eomorphe \`a $\Bbb S^2 \times \Bbb S^1$.
\endproclaim

\medskip

Une surface $F$ plong\'ee dans $M$ est {\it proprement plong\'ee}
si $\partial F$ co\"{\i}ncide avec $F \cap \partial M$ 
ou si $F$ est dans $\partial M$. 
Une surface $F$ proprement plong\'ee dans $M$ est 
soit {\it bilat\`ere} soit {\it unilat\`ere}~;
en d'autres termes, $F$  admet un voisinage r\'egulier dans $M$ 
qui est hom\'eomorphe respectivement 
soit \`a un $I$-fibr\'e trivial soit \`a un $I$-fibr\'e non trivial, 
de base $F$.
Si $F$ est dans $\partial M$, nous convenons que $F$ est bilat\`ere
(ce qui est coh\'erent avec le corollaire 1.10 de \cite{Hemp--76}). 
Lorsque $M$ est orientable, alors $F$ est bilat\`ere 
si et seulement si $F$ est orientable 
(voir \cite{SeTh--34}, \S \ 76, th\'eor\`eme III).

  Une surface proprement plong\'ee $F$ est par d\'efinition  {\it compressible} 
dans les situations suivantes~:
\roster
\item"(i)"
$F$ borde une $3$-boule d'homotopie plong\'ee dans $M$, 
\item"(ii)"
$F$ est une $2$-cellule dans $\partial M$,
\item"(iii)" 
$F$ est une $2$-cellule dans $M$ 
et il existe une $3$-boule d'homotopie $X$ dans $M$ 
telle que $\partial X \subset F \cup \partial M$, 
\item"(iv)"
il existe une $2$-cellule  $D \subset M$ telle que $D \cap F = \partial D$ 
et telle que $\partial D$ n'est pas contractible dans $F$~; 
\endroster
et $F$ est {\it incompressible} sinon.

   Soit $F$ une surface, distincte de $\Bbb D^2$ et $\Bbb S^2$,
proprement plong\'ee dans $M$.
Si l'homo\-morphisme $\pi_1(F) \longrightarrow \pi_1(M)$ est injectif,
alors la surface $F$ est incompressible 
(voir la situation (iv)).
Si $F$ est bilat\`ere, la r\'eciproque est vraie~;
c'est une application du th\'eor\`eme du lacet
\footnote{
En anglais~: \lq\lq loop theorem\rq\rq , ou encore \lq\lq Dehn+loop theorem\rq\rq .
Voir le corollaire 6.2 de \cite{Hemp--76}
}
\hskip-.1cm
.

\medskip

\proclaim{Formulation de Stallings du th\'eor\`eme de Grushko}
Soit $M$ une $3$-vari\'et\'e telle que
toute composante connexe de $\partial M$ est
incompressible. S'il existe une d\'ecomposition non triviale en
produit libre $\pi_1(M) = \Gamma_1 \ast \Gamma_2$, alors il existe
une d\'ecomposition en somme connexe $M = M_1 \sharp\, M_2$ telle
que $\Gamma_1 = \pi_1(M_1)$ et $\Gamma_2 = \pi_1(M_2)$.
\endproclaim

\demo{R\'ef\'erence et remarque}
Voir le th\'eor\`eme 7.1 dans \cite{Hemp--76}.
Exemple montrant que la condition sur $\partial M$ est n\'ecessaire~: 
un corps \`a anses de genre $g \ge 2$ est une
$3$-vari\'et\'e irr\'eductible, \`a bord compressible, 
dont le groupe fondamental est un produit libre non trivial.
\enddemo

\medskip


Une $3$-vari\'et\'e est {\it $\Bbb P^2$-irr\'eductible} si
elle est irr\'eductible et si elle ne contient aucune surface
plong\'ee bilat\`ere hom\'eomorphe au plan projectif $\Bbb P^2$.
En particulier, lorsque $M$ est orientable, 
$M$ est $\Bbb P^2$-irr\'eductible si et seulement si $M$ est irr\'eductible.

\bigskip
\head{\bf 
III.~Les vari\'et\'es de Seifert
}
\endhead
\medskip

Sur une $3$-vari\'et\'e $M$, 
une {\it structure de Seifert}, ou {\it fibration de Seifert}, 
est un feuilletage en cercles
tel que chaque feuille poss\`ede un voisinage hom\'eomorphe 
\`a un tore plein ou \`a une bouteille de Klein pleine
feuillet\'e \lq\lq de mani\`ere standard
\footnote{
\lq\lq Standard\rq\rq \ se r\'ef\`ere \`a une condition de r\'egularit\'e
qui est automatiquement v\'erifi\'ee,
par un th\'eor\`eme d\^u \`a Epstein \cite{Epst--72}~;
voir aussi le \S~3 de \cite{Sco--83b}.
Les bouteilles de Klein n'apparaissent \'evidemment pas
lorsque $M$ est orientable.
}
\rq\rq .
Les {\it fibres r\'eguli\`eres} de la fibration de Seifert
sont les feuilles pour lesquelles l'application d'holonomie 
(ou application de premier retour sur un petit disque transverse \`a la
feuille) est l'identit\'e. 

Une {\it vari\'et\'e de Seifert} est une $3$-vari\'et\'e 
qui poss\`ede une structure de Seifert.

\medskip

Rappelons  d'abord que le groupe fondamental d'une
vari\'et\'e de Seifert non simplement connexe n'est jamais cci.
C'est une cons\'equence \'el\'ementaire 
du fait que la \lq\lq structure fibr\'ee passe de $M$ \`a $\Gamma$\rq\rq .

\medskip

\proclaim{1.\ Lemme}
Soient $M$ une $3$-vari\'et\'e sur laquelle il existe une structure
de Seifert et $X$ l'espace des orbites. Alors il existe une suite exacte courte
$$ 
1 \quad \longrightarrow \quad 
K \quad \longrightarrow \quad
\pi_1(M) \quad \longrightarrow \quad
\pi_1^{\text{orb}}(X) \quad \longrightarrow \quad 
1
$$
o\`u $K$ est le sous-groupe de $\pi_1(M)$ engendr\'e par la classe d'homotopie
d'une fibre r\'eguli\`ere et o\`u $\pi_1^{\text{orb}}(X)$ d\'esigne 
le groupe fondamental de $X$ au sens des orbi\'et\'es.
De plus, si le rev\^etement universel de $M$ 
n'est pas hom\'eomorphe \`a $\Bbb S^3$,
alors $K$ est un groupe cyclique infini.
\endproclaim

\demo\nofrills{Pour la d\'emonstration~: \usualspace} 
voir le lemme 3.2 de \cite{Sco--83b}.
\hfill $\square$
\enddemo

\medskip

\proclaim{2.\ Proposition} Le groupe fondamental d'une vari\'et\'e
de Seifert non simplement connexe poss\`ede une classe de
conjugaison finie autre que $1$.
\endproclaim

\demo{D\'emonstration} Avec les notations du lemme pr\'ec\'edent,
ou bien $\pi_1(M)$ est fini non r\'eduit \`a $1$, et il n'y a rien \`a montrer,
ou bien le groupe $K$ est cyclique infini, et en particulier
contient des classes de conjugaison finies autres que $1$.
\hfill $\square$
\enddemo

Dans la suite, nous montrons dans quelle mesure, r\'eciproquement, 
une $3$-vari\'et\'e \`a groupe fondamental infini non cci 
est essentiellement une vari\'et\'e de Seifert. 
Apr\`es quelques rappels et des pr\'eliminaires sur les groupes cci,
nous traitons successivement les vari\'et\'es orientables irr\'eductibles,
orientables, non orientables $\Bbb P^2$-irr\'eductibles, et non orientables.

\bigskip
\head{\bf 
IV.~Rappels sur
la conjecture des fibr\'es de Seifert (cas orientable)
\\
et sur un r\'esultat de Hempel-Jaco
}
\endhead
\medskip 

Nos preuves utilisent deux ingr\'edients cruciaux. 
D'abord la \lq\lq conjecture des fibr\'es  de Seifert\rq\rq ,
qui est devenue un th\'eor\`eme
gr\^ace aux travaux de 
Waldhausen, Gordon-Heil, Jaco-Shalen, Scott, Mess, Tukia, Gabai, Casson-Jungreis
(pour l'histoire, voir \cite{Gaba--92}, pages 507--508),
ainsi que Whitten et Heil (pour le cas non orientable).

\medskip

\proclaim{Th\'eor\`eme (conjecture des fibr\'es  de Seifert, cas orientable)} 
Soit $M$ une $3$-vari\'et\'e orientable irr\'e\-ductible. 
Si le groupe fondamental de $M$ 
contient un sous-groupe normal cyclique infini, 
alors $M$ est une vari\'et\'e de Seifert.
\endproclaim

\medskip

\demo{Remarques} (i) Soit $M$ une $3$-vari\'et\'e obtenue 
par supression d'une $3$-boule ouverte 
dans une vari\'et\'e de Seifert $\hat M$.
D'une part, le bord de $M$ poss\`ede une composante connexe $\Bbb S^2$,
ce qui exclut que $M$ soit une vari\'et\'e de Seifert~;
d'autre part, $\pi_1(M)$ est isomorphe \`a $\pi_1(\hat M)$. 
On ne peut donc pas supprimer l'hypoth\`ese d'irr\'eductibilit\'e dans le
th\'eor\`eme pr\'ec\'edent.

(ii) Pour le cas des vari\'et\'es non orientables, voir le \S~VIII ci-dessous.
\enddemo

Voici ensuite ce que nous voulons rappeler ici 
d'un th\'eor\`eme de Hempel et Jaco. 
Pour l'\'enonc\'e complet, 
voir \cite{HeJa--72}, ou le th\'eor\`eme 11.1 de \cite{Hemp--76}  
(dans lequel il faut prendre garde
\`a l'hypoth\`ese implicite $N \ne 1$).

\medskip

\proclaim{Th\'eor\`eme (Hempel-Jaco)}
Soit $M$ une $3$-vari\'et\'e $\Bbb P^2$-irr\'eductible
dont le groupe fondamental $\Gamma$ s'ins\`ere dans une suite exacte courte
$$
1 \quad \longrightarrow \quad 
K \quad \longrightarrow \quad
\Gamma \quad \longrightarrow \quad
Q \quad \longrightarrow \quad 
1
$$
avec $K  \ne 1$, $K$ de type fini et $Q$ infini. 
Si $K$ n'est pas isomorphe \`a $\Bbb Z$,
alors l'une au moins des assertions suivantes est vraie
\roster
\item"(i)"    $M$ est un fibr\'e au dessus du cercle 
   \`a fibre une surface $F$~;
\item"(ii)"     il existe un $I$-fibr\'e non trivial $E$
   sur une surface $F$ tel que $M$ soit hom\'eo\-morphe \`a la
   r\'eunion de deux copies de $E$ recoll\'ees sur leur bord commun.
\endroster
Dans les deux cas, 
$K$ est isomorphe \`a un sous-groupe d'indice fini dans
$\pi_1(F)$~; en particulier, $K$ est un groupe de surface.
\endproclaim

\bigskip
\head{\bf 
V.~Groupes cci
}\endhead
\medskip

Soit $\Gamma = \Gamma_1 \ast \Gamma_2$ un produit libre de deux groupes
dont aucun n'est  r\'eduit \`a~$1$.
Si $\Gamma_1$ et $\Gamma_2$ sont d'ordre $2$,
alors $\Gamma$ est un groupe di\'edral infini et n'est donc pas cci.

\medskip

\proclaim{3.\ Proposition}
(i) Un produit libre $\Gamma = \Gamma_1 \ast \Gamma_2$ 
tel que $\abs{\Gamma_1} \ge 3$ et $\abs{\Gamma_2} \ge 2$
est un groupe~cci.

Soit $\Gamma = \Gamma_1 \ast_{\Gamma_0} \Gamma_2$
un produit libre avec amalgamation relativement \`a des inclusions
$\Gamma_0 \subset \Gamma_1$ et $\Gamma_0 \subset \Gamma_2$
d'indices $[\Gamma_1 : \Gamma_0] \ge 3$ et $[\Gamma_2 : \Gamma_0] \ge 2$.

(ii)  Si l'un au moins des groupes $\Gamma_1$, $\Gamma_2$ est cci,
le groupe $\Gamma$ est aussi cci.

(iii) On suppose que tout sous-groupe de type fini de $\Gamma_0$ 
non r\'eduit \`a un \'el\'ement 
qui est normal \`a la fois dans $\Gamma_1$ et dans $\Gamma_2$ 
est cci.
Alors le groupe $\Gamma$ est aussi cci.

Soient $\Gamma_1$ un groupe, 
$\Gamma_0, \Gamma'_0$ deux sous-groupes de $\Gamma_1$, 
l'un au moins \'etant propre,
$\varphi$ un isomorphisme de $\Gamma_0$ sur $\Gamma'_0$,
et $\Gamma$ l'extension HNN correspondante.

(iv) Si $\Gamma_1$ est cci, le groupe $\Gamma$ est aussi cci.

(v) On suppose que tout sous-groupe de type fini de $\Gamma_0$ 
non r\'eduit \`a un \'el\'ement
qui est normal dans $\Gamma$ 
est cci. 
Alors le groupe $\Gamma$ est aussi cci.

\endproclaim

\demo{D\'emonstration}
Ce sont des cons\'equences des th\'eor\`emes usuels de formes normales.
A~titre d'exemple, 
d\'etaillons l'argument pour (v) 
lorsque $\Gamma_0$ est un sous-groupe propre de $\Gamma_1$.
Notons $t$ la lettre stable de l'extension HNN,
et $C_{\Delta}(\delta)$ la classe de conjugaison dans un groupe $\Delta$
d'un \'el\'ement $\delta \in \Delta$.
Soit $\gamma \in \Gamma$, $\gamma \ne 1$~;
il s'agit de montrer que la classe $C_{\Gamma}(\gamma)$ est infinie.

Si $\gamma \notin \Gamma_1$, la classe $C_{\Gamma}(\gamma)$
est infinie par le lemme 2.1 de \cite{Stal}.
On suppose d\'esormais que $\gamma \in \Gamma_1$.
Si $\gamma \notin \Gamma_0$, les \'el\'ements $t^{-n} \gamma t^n$
(o\`u $n = 1,2,3,\hdots$) sont distincts deux \`a deux, 
car ce sont des \'ecritures de formes normales.
On suppose d\'esormais que $\gamma \in \Gamma_0$.

Si la classe $C_{\Gamma}(\gamma)$ \'etait finie,
elle serait contenue dans $\Gamma_0$ 
vu les cas d\'ej\`a trait\'es. 
Le groupe engendr\'e par $C_{\Gamma}(\gamma)$ serait donc 
un sous-groupe de type fini de $\Gamma_0$ qui serait normal dans $\Gamma$, 
et il poss\`ederait un \'el\'ement distinct de $1$ (\`a savoir $\gamma$)
ayant une classe de conjugaison $C_{\Gamma_0}(\gamma)$ finie,
contrairement aux hypoth\`eses. 
Par suite la classe $C_{\Gamma}(\gamma)$ est infinie.
\hfill $\square$
\enddemo

\demo{Exemples} Soit $F$ une surface compacte connexe qui n'est pas
hom\'eomorphe \`a l'une des surfaces de la liste
\par
\centerline{
$2$-disque, anneau, ruban de M\"obius,
$\Bbb S^2$, $\Bbb P^2$, $2$-tore, bouteille de Klein, 
}
\par\noindent
de sorte que le groupe fondamental $\pi_1(F)$ 
n'est ni fini ni virtuellement ab\'elien.
Si le bord de $F$ n'est pas vide, alors $\pi_1(F)$ est libre non ab\'elien,
donc cci par l'assertion (i) de la proposition pr\'ec\'edente.
Si le bord de $F$ est vide, alors $F$ est une surface close,
orientable de genre au moins $2$ ou non orientable de genre au moins $3$,
donc $\pi_1(F)$ est cci par l'assertion~(ii).
\enddemo

\medskip

Il r\'esulte du th\'eor\`eme de Kneser-Milnor que,
\lq\lq en g\'en\'eral\rq\rq , 
le groupe fondamental d'une $3$-vari\'et\'e d\'ecomposable est cci~;
comme bien d'autres \lq\lq \'enonc\'es g\'en\'eraux\rq\rq \
concernant les $3$-vari\'et\'es, 
celui-ci est contredit par quelques exceptions
qui contribuent \`a allonger les arguments qui suivent.

Pour la commodit\'e du lecteur, 
nous collectons ici trois lemmes concernant les groupes cci.

\medskip

\proclaim{4.\ Lemme} Dans un groupe cci, 
tout sous-groupe d'indice fini est cci.
\endproclaim

\demo{D\'emonstration} Soient $G$ un groupe cci et 
$H$ un sous-groupe d'indice fini.
Choisissons un sous-ensemble fini $T$ de $G$ tel que 
$G$ soit la r\'eunion disjointe $\sqcup_{t \in T} tH$.
Soit  $h \in H$ un \'el\'ement 
dont la classe de conjugaison $C_H(h) = \{h_1,\hdots,h_k\}$ dans $H$ 
est finie~;
les \'el\'ements $th_it^{-1}$ ($t \in T, 1 \le i \le k$)
constituent une \'enum\'eration (peut-\^etre avec r\'ep\'etitions)
de sa classe $C_G(h)$ dans $G$~;
en particulier, celle-ci est finie, donc $h=1$.
\hfill $\square$
\enddemo

La r\'eciproque du lemme 4 n'est pas correcte, 
comme on le voit en consid\'erant
le produit direct d'un groupe cci et d'un groupe fini non r\'eduit \`a un
\'el\'ement.   

\medskip

\proclaim{5.\ Lemme} 
Un groupe sans torsion qui poss\`ede un sous-groupe d'indice fini cci
est lui-m\^eme cci.
\endproclaim

\demo{D\'emonstration} 
Soit $G$ un groupe poss\'edant un sous-groupe d'indice fini $H$ qui est cci.
Il suffit de montrer que, si $G$ n'est pas cci, alors $G$ poss\`ede de la
torsion.

Soit donc $g \in G$, $g \ne 1$, un \'el\'ement dont la classe de conjugaison
est finie.  Il existe un entier $n > 1$ tel que $g^n \in H$, 
car $H$ est d'indice fini. De plus $g^n = 1$, car $H$ est cci.
\hfill $\square$
\enddemo
 
\medskip

\proclaim{6.\ Lemme} Soient $G$ un groupe et $H$ un sous-groupe d'indice
$2$~; on suppose que $H$ est cci.
Alors ou bien $G$ est cci, 
ou bien $G$ est  produit direct de $H$ et d'un groupe d'ordre $2$.
\endproclaim

\demo {D\'emonstration} Notons $F_G$ la r\'eunion des classes
de conjugaison finies de $G$~; c'est un sous-groupe distingu\'e de $G$.
De plus $F_G$ est d'ordre au plus $2$, 
car $F_G \cap H \subset F_H = 1$.
Si $F_G$ est r\'eduit \`a $1$, alors $G$ est cci.

Supposons d\'esormais que $F_G$ soit d'ordre $2$,
et notons $f$ son \'el\'ement distinct de l'identit\'e.
Alors $f$ est d'ordre $2$ et central dans $G$,
de sorte que $G$ est produit direct de $H$ et $F_G$. 
\hfill $\square$
\enddemo

\bigskip
\head{\bf 
VI.~Le th\'eor\`eme principal pour les vari\'et\'es orientables irr\'eductibles
}\endhead
\medskip

Avec une $3$-vari\'et\'e $M$, il convient de consid\'erer 
la vari\'et\'e $\hat M$  obtenue \`a partir de $M$ 
en lui recollant une $3$-cellule  
le long de chaque composante de bord hom\'eomorphe \`a $\Bbb S^2$
par un hom\'eomorphisme renversant l'orientation. 
C'est une cons\'equence imm\'ediate du th\'eor\`eme de Seifert-Van Kampen 
que l'inclusion $M \subset \hat M$ induit un isomorphisme 
$\pi_1(M) \approx \pi_1(\hat M)$.

Soit $M$ une vari\'et\'e et
$$
M \, = \, M_1 \sharp\, \cdots \sharp\, M_k \sharp\, 
          B_1 \sharp\, \cdots \sharp\, B_l \sharp\, 
          C_1 \sharp\, \cdots \sharp\, C_m
$$
une d\'ecomposition de Kneser-Milnor de $M$,
avec les $M_i$ ind\'ecomposables et non simplement connexes,
les $B_j$ des $3$-boules
et les $C_k$ des $3$-sph\`eres non standard
(c'est-\`a-dire 
\footnote{
Nous ne supposons rien ici concernant
la conjecture de Poincar\'e.
}
des $3$-sph\`eres d'homotopie non hom\'eomorphes \`a $\Bbb S^3$).
Notons que 
$$
\hat M \, = \, \left\{
\aligned
M_1 \sharp\, \cdots \sharp\, M_k \sharp\, 
          C_1 \sharp\, \cdots \sharp\, C_m 
          \qquad &\text{si} \quad k+m \ge 1,\\
\Bbb S^3 \hskip2.5cm &\text{si} \quad k+m = 0 .
\endaligned
\right.
$$
La {\it vari\'et\'e de Poincar\'e} associ\'ee \`a $M$ 
est la vari\'et\'e d\'efinie par
$$
\Cal P (M) \, = \, \left\{
\aligned 
M_1 \sharp\, \cdots \sharp\, M_{k} \qquad &\text{si} \quad k \ge 1, \\
\Bbb S^3 \hskip1.5cm &\text{si} \quad k = 0 .
\endaligned
\right.
$$
Il r\'esulte \`a nouveau du th\'eor\`eme de Seifert-Van Kampen que
les groupes fondamentaux de $M$ et $\Cal P (M)$ sont isomorphes.
De plus, toute sous-vari\'et\'e compacte contractile 
de dimension trois de $\Cal P (M)$ 
est hom\'eomorphe \`a un $3$-disque~; 
voir le d\'ebut du chapitre 10 de \cite{Hemp--76}. 

Pour le lecteur qui accepterait la conjecture de Poincar\'e,
il conviendrait de lire $\hat M$ au lieu de $\Cal P (M)$ ci-dessous.

\medskip

Les quatre lemmes qui suivent pr\'eparent la preuve du th\'eor\`eme 12.

\medskip

\proclaim{7.\ Lemme} 
   (i) Soit $M$  une vari\'et\'e irr\'eductible 
qui n'est pas une $3$-boule. Alors $\hat M = M$~;
autrement dit, aucune composante connexe du bord $\partial M$ 
n'est une $2$-sph\`ere.

   (ii) Si $M$ est une vari\'et\'e irr\'eductible non simplement connexe, 
alors $\Cal P (M) = M$.
\endproclaim

\demo{D\'emonstration} L'assertion (i) provient du fait que
remplacer une vari\'et\'e $M$ par sa somme connexe avec une $3$-boule
revient \`a cr\'eer dans $M$ une composante de bord $\Bbb S^2$.
L'assertion (ii) est une cons\'equence imm\'ediate des d\'efinitions.
\hfill $\square$
\enddemo

\medskip

\proclaim{8.\ Lemme} Soient $G$ un groupe et 
$H$ un sous-groupe d'indice fini de $G$ ayant un centre infini. 
Alors $H$ contient un sous-groupe distingu\'e d'indice fini de $G$ 
ayant un centre infini.
\endproclaim

\demo{D\'emonstration} 
Soit $N$ le noyau de l'homomorphisme naturel $\alpha$ de $G$ 
dans le groupe des permutations du $G$-ensemble fini $G/H$. 
Alors $N$ est normal et d'indice fini dans $G$. 
De plus, la restriction de $\alpha$ au centre de $H$ 
a un noyau qui est infini (car d'indice fini dans le centre de $H$)
et qui est contenu dans le centre de $N$. 
\hfill $\square$
\enddemo

\medskip

\proclaim{9.\ Lemme} Soient $G$ un groupe contenant un \'el\'ement $g$
d'ordre infini \`a centralisateur $Z_G(g)$  d'indice fini dans $G$. 
Alors $G$ contient un sous-groupe ab\'elien de type fini qui est infini et
distingu\'e dans $G$.
\endproclaim

\demo{D\'emonstration} 
Par le lemme pr\'ec\'edent, 
$G$ contient un sous-groupe distingu\'e d'indice fini $N$ 
\`a centre $Z(N)$ infini~;
soient $t_1, \hdots,t_n$ des repr\'esentants de $G$ modulo $N$.
Soient $z \ne e$ un \'el\'ement de $Z(N)$ 
et $K$ la cl\^oture normale de $z$ dans $G$.

Le groupe $K$ est engendr\'e par 
$t_1 z t_1^{-1}, \hdots, t_n z t_n^{-1}$,
de sorte que $K$ est de type fini.
Le groupe $Z(N)$ est caract\'eristique dans $N$, qui est normal dans $G$,
de sorte que $Z(N)$ est normal dans $G$~; il en r\'esulte que $K$ est
contenu dans $Z(N)$, et en particulier que $K$ est ab\'elien.
Enfin $K$ est normal dans $G$, par d\'efinition.
\hfill $\square$
\enddemo

\medskip

\proclaim{10.\ Lemme} 
Soit $N$ l'espace total d'un fibr\'e en $2$-tores sur le cercle. 
Si $N$ est orientable et si le groupe $\pi_1(N)$  n'est pas cci,
alors $N$ est une vari\'et\'e de Seifert.
\endproclaim

\demo{D\'emonstration} 
Le groupe fondamental $\pi_1(N)$ est produit semi-direct 
de $\Bbb Z^2$ par $\Bbb Z$ relativement \`a un automorphisme
$\phi \in \operatorname{Aut}(\Bbb Z^2) = GL(2,\Bbb Z)$.
Plus pr\'ecis\'ement, $\phi$ est induit sur le groupe fondamental
par un hom\'eomorphisme du $2$-tore qui pr\'eserve l'orientation
(car $N$ est orientable), de sorte que $\phi \in SL(2,\Bbb Z)$.
Il y a {\it a priori} trois cas \`a consid\'erer selon que $\phi$ est
elliptique (c'est-\`a-dire \`a valeurs propres non r\'eelles),
parabolique (\`a valeurs propres $1$ ou $-1$),
ou hyperbolique 
(\`a valeurs propres r\'eelles de modules diff\'erents de $1$).

{\it Si $\phi$ est elliptique,} $\phi$ est d'ordre fini, 
de sorte que $N$ est un fibr\'e de Seifert
(lemme II.5.4 de \cite{JaSh--79}).

{\it Si $\phi$ est parabolique,} alors $\phi$ est conjugu\'e  
\`a une matrice de la forme 
$\pm \left( \matrix 1 & k \\ 0 & 1 \endmatrix \right)$.
Sur le tore, $\phi$ d\'efinit soit un twist de Dehn (signe $+$)
soit un twist de Dehn compos\'e avec la sym\'etrie centrale (signe $-$),
donc laisse invariant un feuilletage en cercles. 
Par suite la vari\'et\'e $N$ admet un feuilletage en cercles~;
c'est donc une vari\'et\'e de Seifert.

{\it Si $\phi$ est hyperbolique,} 
$\phi$ n'a pas de vecteur propre dans $\Bbb Z^ 2$~;
il est alors facile de v\'erifier que
le groupe $\pi_1(N) = \Bbb Z^2 \rtimes_{\phi} \Bbb Z$ est cci,
et ce cas n'entre donc pas dans les hypoth\`eses du lemme.
\hfill $\square$
\enddemo

\medskip  

Rappelons que, \`a une extension de groupes
$$
1 \quad \longrightarrow \quad
A  \quad \longrightarrow \quad
B  \quad \overset{\pi}\to{\longrightarrow} \quad
C  \quad \longrightarrow \quad 1
$$
avec $A$ ab\'elien,
on associe naturellement  l'homomorphisme
$\theta : C \longrightarrow \operatorname{Aut}(A)$ d\'efini comme suit~;
pour $c \in C$, l'automorphisme $\theta(c)$
est la restriction \`a $A$ de la conjugaison par 
un \'el\'ement de $\pi^{-1}(c)$~;
voir par exemple \cite{Rotm--95}, pages 178 et suivantes.
Il r\'esulte de cette d\'efinition que  le noyau de $\theta$
est l'image par $\pi$ du centralisateur de $A$ dans~$B$.

\medskip

\proclaim{11.\ Lemme}
Soit $\Gamma$ un $3$-groupe sans torsion tel qu'il existe une extension
$$
1 \quad \longrightarrow \quad
K  \quad \longrightarrow \quad
\Gamma  \quad \overset{\pi}\to{\longrightarrow} \quad
Q  \quad \longrightarrow \quad 1
$$
avec $Q$ fini, 
$K$ ab\'elien libre de rang trois,
et $K$ maximal parmi les sous-groupes ab\'eliens libres
de type fini distingu\'es dans $\Gamma$.  
Alors l'homomorphisme $\theta : Q \longrightarrow GL(3,\Bbb Z)$ 
associ\'e \`a l'extension est injectif.
\endproclaim

\demo{D\'emonstration} 
Le noyau $Q_0$  de $\theta$ et son image inverse $\Gamma_0 = \pi^{-1}(Q_0)$
s'ins\`erent dans une suite exacte
$$
1 \quad \longrightarrow \quad
K  \quad \longrightarrow \quad
\Gamma_0  \quad \overset{\pi_0}\to{\longrightarrow} \quad
Q_0  \quad \longrightarrow \quad 1
$$
dont l'homomorphisme associ\'e $\theta_0$ 
est l'homomorphisme constant de $Q_0$ dans $GL(3,\Bbb Z)$~;
il en r\'esulte que $K$ est central dans $\Gamma_0$.
Le th\'eor\`eme 12.10 de \cite{Hemp--76} implique que
$\Gamma_0$ est isomorphe \`a $\Bbb Z^3$,
donc que $K = \Gamma_0$ par hypoth\`ese de maximalit\'e,
de sorte que le quotient $Q_0$ est r\'eduit \`a $1$.
\hfill $\square$
\enddemo

\medskip

\proclaim{12.\ Th\'eor\`eme} 
Soit $M$ une $3$-vari\'et\'e orientable irr\'eductible  
dont le groupe fondamental $\Gamma$ est infini. 
Alors $\Gamma$ n'est pas cci si et seulement si 
$M$ est une vari\'et\'e de Seifert.
\endproclaim

\demo{D\'emonstration} Vu la proposition 2, nous pouvons supposer
que $\Gamma$ n'est pas cci, c'est-\`a-dire qu'il contient
un \'el\'ement $\gamma \ne 1$ dont le centralisateur 
$Z_{\Gamma}(\gamma)$ est d'indice fini dans $\Gamma$~;
il s'agit de montrer que $M$ est une vari\'et\'e de Seifert.

Le groupe $\Gamma$ est sans torsion 
car $M$ est un espace d'Eilenberg-MacLane. 
(En effet, $\pi_2(M)$ est trivial par le th\'eor\`eme de la sph\`ere~;
d\'etails au corollaire 9.9 de \cite{Hemp--76}.) 
Le lemme 9 montre qu'il existe un sous-groupe normal $K$ dans $\Gamma$ 
qui est ab\'elien, de type fini et infini~;
d'o\`u une suite exacte courte
$$
1 \quad \longrightarrow \quad
K  \quad \longrightarrow \quad
\Gamma  \quad \longrightarrow \quad
Q  \quad \longrightarrow \quad 1
$$
(avec $Q = \Gamma / K$).
La liste est connue des groupes ab\'eliens de type fini
qui peuvent \^etre sous-groupes de $3$-groupes 
(th\'eor\`eme 9.13 de \cite{Hemp--76}).
Comme $K$ est de plus sans torsion, il est isomorphe \`a l'un des groupes
$\Bbb Z$, $\Bbb Z \oplus \Bbb Z$, 
et $\Bbb Z \oplus \Bbb Z \oplus \Bbb Z$.

\smallskip
\noindent {\it Premier cas : $K$ est cyclique infini.}

Il r\'esulte de la  \lq\lq conjecture\rq\rq \ des fibr\'es de Seifert 
(voir plus haut)  que $M$ est une vari\'et\'e de Seifert.

\smallskip
\noindent {\it Deuxi\`eme cas : 
$K \approx \Bbb Z \oplus \Bbb Z$ et $Q$ est fini.} 

La vari\'et\'e $\Cal P (M)$ est un $I$-fibr\'e sur une surface close
par le th\'eor\`eme 10.6 de \cite{Hemp--76}, qui ne peut \^etre qu'un
$2$-tore ou une bouteille de Klein 
(car de groupe fondamental virtuellement $\Bbb Z^2$)~;
de plus $\Cal P (M) = M$ (lemme 7).
Il en r\'esulte que $M$ est une vari\'et\'e de Seifert.

\smallskip
\noindent {\it Troisi\`eme cas : 
$K \approx \Bbb Z \oplus \Bbb Z$ et $Q$ est infini.} 

Il r\'esulte du th\'eor\`eme de Hempel-Jaco cit\'e plus haut que 
l'une au moins des assertions suivantes est vraie~:
\roster
\item"(i)" $M$ est un fibr\'e  sur le cercle dont la fibre $F$
   est un tore~; 
\item"(ii)" $M$ est un recollement convenable le long d'une surface $F$
   de deux $I$-fibr\'es.
\endroster
En particulier, $M$ est une vari\'et\'e de Haken
(la surface $F$ du pr\'esent contexte est bilat\`ere incompressible)
et $M$ poss\`ede un rev\^etement fini~$N$
qui est l'espace total d'un fibr\'e en tores sur le cercle.
Le lemme 6 implique que le groupe infini $\pi_1(N)$ n'est pas cci,
et le lemme 10 que $N$ est une vari\'et\'e de Seifert.

La vari\'et\'e $M$ est de Seifert, car elle est de Haken 
et elle poss\`ede un rev\^etement fini qui est une vari\'et\'e de Seifert
(th\'eor\`eme II.6.3 de \cite{JaSh--79}).

\smallskip
\noindent {\it Quatri\`eme cas : 
$K \approx \Bbb Z^3$.} 

Le groupe $Q$ est fini par le th\'eor\`eme de Hempel-Jaco
(et le fait que $\Bbb Z^3$ n'est pas un groupe de surface !).
Quitte \`a remplacer $K$ par un sur-groupe d'indice fini,
on peut supposer $K$ maximal parmi les sous-groupes de $\Gamma$ qui sont 
distingu\'es, ab\'eliens et de type fini.
L'homomorphisme
$$
\theta \, : \, Q \longrightarrow \operatorname{Aut}(\Bbb Z^3)
\approx GL_3(\Bbb Z) 
$$
associ\'e \`a l'extension 
$1 \longrightarrow K = \Bbb Z^3 \longrightarrow
\Gamma \longrightarrow Q \longrightarrow 1$
est injectif (lemme 11).
Un argument qu'on trouve dans \cite{Sco--83b} (page 444)
montre alors que $\Gamma$ poss\`ede un sous-groupe normal cyclique infini.
Il en r\'esulte que $M$ est une vari\'et\'e de Seifert
(\lq\lq conjecture\rq\rq \ des fibr\'es de Seifert).
\hfill $\square$
\enddemo

\bigskip
\head{\bf 
VII.~Le th\'eor\`eme principal pour les vari\'et\'es orientables 
}\endhead
\medskip

Nous pouvons maintenant traiter le cas d'une 3-vari\'et\'e
orientable en toute g\'en\'eralit\'e.
Le prix qu'implique la consid\'eration  des vari\'et\'es r\'eductibles
est la pr\'esence a priori possible de $3$-sph\`eres non standard
dans la d\'ecomposition de Kneser-Milnor de $M$,
d'o\`u l'apparition de la vari\'et\'e de Poincar\'e $\Cal P(M)$.

\medskip

\proclaim{13.\ Th\'eor\`eme} Soit $M$ une $3$-vari\'et\'e orientable 
et $\Gamma$ son groupe fondamental~;
on suppose $\Gamma$ infini et non di\'edral infini. 
Alors $\Gamma$ n'est pas cci si et seulement si $\Cal P (M)$ est une vari\'et\'e
de Seifert.

En particulier, pour le groupe fondamental infini $\Gamma$
d'une $3$-vari\'et\'e orientable,
les trois propri\'et\'es suivantes sont \'equivalentes~:
\roster
\item"(i)" $\Gamma$ est groupe fondamental d'une vari\'et\'e de Seifert~; 
\item"(ii)" $\Gamma$ poss\`ede un sous-groupe normal cyclique infini~;
\item"(iii)" $\Gamma$ n'est pas cci.
\endroster
\endproclaim

\demo{Remarques} (i) Au th\'eor\`eme 13,  il n'est pas possible de remplacer 
$\Cal P (M)$ par $M$. 
En effet, soient $M_1$ une vari\'et\'e de Seifert orientable
\`a groupe fondamental $\Gamma$ infini et $C$ une $3$-sph\`ere non standard~;
posons $M = M_1 \sharp\,  C$, de sorte que $\Cal P (M) = M_1$.
Alors  $\Gamma = \pi_1(M_1) = \pi_1(M)$ n'est pas cci (proposition 2), 
$M$ n'est pas une vari\'et\'e de Seifert \cite{Sco--83a}
et $\Cal P (M)$ en est une.

(ii) Acceptons pour cette remarque la conjecture de Poincar\'e.
Soit $M$ une $3$-vari\'et\'e orientable de la forme
$M_1 \sharp\, M_2$, o\`u les vari\'et\'es $M_1,M_2$ sont premi\`eres
et les groupes $\pi_1(M_1), \pi_1(M_2)$ d'ordre $2$.
Alors $M$ est n\'ecessairement $\Bbb P^3 \sharp\, \Bbb P^3$,
qui admet une fibration de Seifert. On peut donc supprimer l'hypoth\`ese
$\Gamma \not\approx D_{\infty}$ dans la premi\`ere assertion du th\'eor\`eme 13.
\enddemo

\demo{D\'emonstration} Vu la proposition 2, il reste \`a montrer que,
si $\Gamma$ n'est pas cci, 
alors $\Cal P (M)$ est une vari\'et\'e de Seifert.
Nous pouvons supposer $M$ r\'eductible
gr\^ace au th\'eor\`eme 12.
Si $M$ est de plus ind\'ecomposable, alors $M$ est hom\'eomorphe
\`a $\Bbb S^2 \times \Bbb S^1$ et il n'y a plus rien \`a montrer.
Nous pouvons donc supposer $M$ d\'ecomposable, et plus pr\'ecis\'ement que
la d\'ecomposition de Kneser-Milnor de $\hat M$ est de la forme
$$
\hat M \, = \, M_1 \sharp\, \cdots \sharp\, M_p \sharp\,
               C_1 \sharp\, \cdots \sharp\, C_q
\hskip1cm (p > 0)
$$
o\`u les $M_i$ sont des vari\'et\'es orientables irr\'eductibles \`a groupes
fondamentaux non r\'eduits \`a $1$ 
et o\`u les $C_j$ sont des sph\`eres non standard.
Si $p=1$, alors $\Cal P(M) = M_1$ est une vari\'et\'e irr\'eductible, 
et le th\'eor\`eme 12 s'applique, 
de sorte que $\Cal P(M)$ est une vari\'et\'e de Seifert.
Si $p \ge 2$, alors
$$
\Gamma \, = \, \Gamma_1 \ast \cdots \ast \Gamma_p
\hskip1cm (\Gamma_i = \pi_1(M_i))
$$
est un produit libre non banal.
Notons que $p = 2$ 
et que les deux groupes $\Gamma_1$, $\Gamma_2$ sont d'ordre~$2$,
car $\Gamma$ n'est pas cci (proposition 3)~;
donc $\Gamma$ est produit libre de deux groupes d'ordre $2$,
c'est un groupe di\'edral infini.
\hfill $\square$
\enddemo

\bigskip
\head{\bf 
VIII.~Rappel sur la conjecture des fibr\'es de Seifert (cas non orientable) 
}\endhead
\medskip

Notons 
\footnote{
Plus pr\'ecis\'ement, 
on consid\`ere deux copies $P_1,P_2$ de la vari\'et\'e $\Bbb P^2 \times I$, 
ainsi que deux $2$-disques plong\'es dans leurs bords
$D_1 \subset \partial P_1$, $D_2 \subset \partial P_2$,
et on recolle le long de $D_1,D_2$ un cylindre plein $\Bbb D^2 \times I$.
}
$\Bbb P$ la somme connexe sur un disque de deux copies de $\Bbb P^2 \times I$.
Le bord de la vari\'et\'e $\Bbb P$ a trois composantes connexes,
deux hom\'eomorphes au plan projectif 
et une hom\'eomorphe \`a la bouteille de Klein $\Bbb K^2$.
Le groupe fondamental de $\Bbb P$ est 
di\'edral infini.
De plus, $\Bbb K^2$ poss\`ede une fibration en cercles telle que
la classe d'homotopie de chaque fibre 
engendre le sous-groupe d'indice deux de $\pi_1(\Bbb P)$~;
convenons qu'une fibration de la composante connexe $\Bbb K^2$
de $\Bbb P$ est {\it sp\'eciale} si elle poss\`ede cette propri\'et\'e.

Sur une $3$-vari\'et\'e {\it non orientable} $M$, 
on d\'efinit comme Heil et Whitten
une {\it structure de Seifert modulo $\Bbb P$.}
Une telle structure consiste 
en la donn\'ee d'une famille $\Cal W = ( K_1, \hdots, K_l)$ 
(possiblement vide)  
de bouteilles de Klein plong\'ees dans l'int\'erieur de $M$, 
disjointes deux \`a deux,
chacune d'entre elles \'etant bilat\`ere
et munie d'une fibration en cercles.
Les adh\'erences dans $M$
des composantes connexes de $M \setminus \Cal W$
sont d'une part 
des copies $P_1, \hdots, P_m$ de $\Bbb P$ 
et d'autre part
des $3$-vari\'et\'es $M_1, \hdots, M_n$ munies de structures de Seifert~;
chaque $K_i$ est dans le bord d'au moins l'un des $P_j$.
Si $K_i$ borde un $P_j$, sa fibration en cercles 
est sp\'eciale dans $P_j$~;
si $K_i$ borde un $M_k$, sa fibration en cercles 
est induite par la fibration de Seifert de $M_k$.

Sur une $3$-vari\'et\'e non orientable 
munie d'une structure de Seifert modulo $\Bbb P$,
les fibres r\'eguli\`eres (d\'efinies comme au \S~III)
d\'efinissent une classe d'homotopie 
qui engendre un sous-groupe de  $\pi_1(M)$
qui est cyclique, normal, et non r\'eduit \`a un \'el\'ement.
On peut montrer que ce sous-groupe est toujours infini \cite{HeWh--94},
mais nous n'utilisons pas ce fait ci-dessous.

Une {\it vari\'et\'e de Seifert modulo $\Bbb P$}
est une $3$-vari\'et\'e non orientable
qui poss\`ede une structure de Seifert modulo $\Bbb P$.

\medskip

\proclaim{14.\ Remarques} 
  (i) Pour tout $\Bbb P^2$ plong\'e dans une $3$-vari\'et\'e $M$,
le groupe fondamental de $\Bbb P^2$ s'injecte dans celui de $M$.

   (ii) Dans le bord d'une vari\'et\'e de Seifert modulo
$\Bbb P$, le nombre de composantes connexes hom\'eomorphes \`a $\Bbb P^2$
est toujours pair.
Une vari\'et\'e de Seifert modulo $\Bbb P$ est une vari\'et\'e de Seifert
si et seulement si son bord ne contient aucun $\Bbb P^2$.

   (iii) Soit $M$ une vari\'et\'e de Seifert modulo $\Bbb P$. 
Si $\pi_1(M)$ n'a pas d'\'el\'ement d'ordre deux,
alors $M$ est une vari\'et\'e de Seifert.
\endproclaim

\demo\nofrills{}
La remarque (i) est le lemme 5.1 de \cite{Hemp--76},
la remarque (ii) d\'ecoule de la d\'efinition, 
et la troisi\`eme remarque r\'esulte des deux premi\`eres.
\enddemo

\medskip 

   Convenons qu'un {\it $\Bbb P^2 \times I$ non standard}
est une $3$-vari\'et\'e homotope 
au produit d'un plan projectif et d'un intervalle,
mais non hom\'eomorphe \`a ce produit 
(si la conjecture de Poincar\'e est vraie, 
il n'existe pas de $\Bbb P^2 \times I$ non standard).
On d\'efinit de m\^eme un {\it $\Bbb P^2 \times \Bbb S^1$ non standard}.

   Voici l'\'enonc\'e du  th\'eor\`eme 1 de
\cite{HeWh--94}.

\medskip

\proclaim{15.\ Th\'eor\`eme 
(conjecture des fibr\'es de Seifert, cas non orientable)}
Soit $M$ une $3$-vari\'et\'e non orientable irr\'eductible
qui ne contient pas de $\Bbb  P^2 \times I$ non standard.
Le groupe fondamental de $M$ contient un sous-groupe normal cyclique 
non r\'eduit \`a un \'el\'ement
si et seulement si
$M$ est ou bien une vari\'et\'e de Seifert modulo $\Bbb P$
ou bien hom\'eomorphe \`a $\Bbb P^2 \times I$.
\endproclaim

\demo{Remarques} 
(i)    Rappelons qu'une $3$-vari\'et\'e non orientable 
\`a groupe fondamental fini 
est toujours homotope \`a un $\Bbb P^2 \times I$ trou\'e
(th\'eor\`eme  9.6 de \cite{Hemp--76}).

(ii) Un $\Bbb P^2 \times \Bbb S^1$ non standard irr\'eductible
doit contenir un $\Bbb P^2 \times I$ non standard.

En effet. Soient $M$ un $\Bbb P^2 \times \Bbb S^1$ non standard irr\'eductible
et soit $\Gamma = \pi_1(M)$.
Soit $C_2$ l'image dans $\Gamma$ du groupe fondamental de $\Bbb P^2$.
Soit $P$ un $\Bbb P^2$ plong\'e bilat\`ere dans $M$,
dont le groupe fondamental s'identifie \`a $C_2$~; 
un tel $P$ existe en vertu d'un r\'esultat d'Epstein 
(th\'eor\`eme 9.8 de \cite{Hemp--76}). 
Distinguons a priori deux cas selon que $P$ est s\'eparant ou non.

Si $P$ est s\'eparant, notons $M_1,M_2$ les composantes connexes de $M \setminus P$
et $\Gamma_1,\Gamma_2$ leurs groupes fondamentaux.
Alors $\Gamma$ est le produit  amalgam\'e $\Gamma_1 \ast_{C_2} \Gamma_2$.
Le centre de $\Gamma$ est donc contenu dans $C_2$~;
or $\Gamma$ est ab\'elien, donc $\pi_1(M) = C_2$, ce qui est absurde.
Ce cas n'appara\^{\i}t donc pas.

Si $P$ n'est pas s\'eparant, posons $M_1 = M \setminus P$ et notons $\Gamma_1$ son
groupe fondamental.
Alors $\Gamma$ est une extension HNN de base $\Gamma_1$ 
au-dessus de deux sous-groupes d'ordre deux.
Notons $t$ la lettre stable de l'extension HNN.
Si $\Gamma_1$ n'\'etait pas d'ordre deux, il existerait $u \in \Gamma_1$
tel que $tut^{-1}$ soit une \'ecriture r\'eduite,
en particulier tel que $tut^{-1} \ne u$,
et ceci est absurde puisque $\Gamma$ est ab\'elien.
Donc $\Gamma_1$ est d'ordre deux.
Par suite $M_1$ est un $\Bbb P^2 \times I$ non standard
en vertu d'un autre r\'esultat d'Epstein 
(th\'eor\`eme 9.6 de \cite{Hemp--76}). 

(iii) L'\'enonc\'e du th\'eor\`eme 15 se simplifie dans le cas
$\Bbb P^2$-irr\'eductible~:
{\it 
Soit $M$ une $3$-vari\'et\'e non orientable $\Bbb P^2$-irr\'eductible.
Le groupe fondamental de $M$ contient un sous-groupe normal cyclique infini
si et seulement si
$M$ est une vari\'et\'e de Seifert.
}
La remarque~(i) implique que le groupe fondamental d'une telle vari\'et\'e $M$
est n\'ecessairement infini.
\enddemo

\bigskip
\head{\bf 
IX.~Le th\'eor\`eme principal 
pour les vari\'et\'es non orientables $\Bbb P^2$-irr\'eductibles
}\endhead
\medskip

Soient $M$ une vari\'et\'e {\it non orientable}
et $\Gamma$ son groupe fondamental.
Notons $\tilde M$ le rev\^etement d'orientation de $M$
et $\tilde \Gamma$ le sous-groupe d'indice $2$ de $\Gamma$ correspondant~;
le th\'eor\`eme 13 s'applique \`a $\tilde M$ et $\tilde \Gamma$.
Rappelons que les vari\'et\'es $\hat M$ et $\Cal P (M)$
ont \'et\'e d\'efinies peu avant le lemme 7.

\medskip

\proclaim{16.\ Lemme} Soient $M$ une $3$-vari\'et\'e non orientable
\`a groupe fondamental $\Gamma$ infini. 
Alors $\Gamma$ n'est pas cci
si et seulement si l'une au moins des assertion suivantes est vraie~:
\roster 
\item"(i)"
$\Cal P(\tilde{M})$ est une vari\'et\'e de Seifert~;
\item"(ii)" $\hat M$ a le type d'homotopie de $\Bbb P^2\times \Bbb S^1$~;
\item"(iii)" $\tilde \Gamma$ est un groupe di\'edral infini.
\endroster
\endproclaim

\demo{Remarque} 
Dans l'assertion (i), on ne peut remplacer
$\Cal P (\tilde M)$ ni par $M$, ni par $\tilde M$, ni par $\Cal P (M)$.
En effet, soit par exemple
$M=(\Bbb P^2\times I)\sharp\, \Bbb (\Bbb P^2\times I)$~;
son rev\^etement d'orientation est 
$\Bbb S^2\times \Bbb S^1$ moins deux boules.
Donc $\Cal P(\tilde{M})$ est une vari\'et\'e de Seifert, 
mais $M$  n'en est pas une, ni $\tilde{M}$~;
de plus $M=\Cal P(M)$.
Il en est de m\^eme pour la vari\'et\'e $\Bbb P$ du \S~VIII,
dont le rev\^etement d'orientation est
un tore plein priv\'e de deux boules.
\enddemo

\demo{D\'emonstration} 
Supposons d'abord que $\Gamma$ n'est pas cci.
Si $\tilde \Gamma$ n'est pas cci, 
$\Cal P(\tilde M)$ est une vari\'et\'e de Seifert 
ou $\tilde \Gamma$ est di\'edral infini par le th\'eor\`eme 13.
On peut donc supposer $\tilde \Gamma$ cci, 
donc $\Gamma$ produit direct de $\tilde \Gamma$ 
et d'un groupe d'ordre $2$ par le lemme 6.

Rappelons qu'un $3$-groupe infini 
poss\`ede des \'el\'ements d'ordre infini.
En effet, 
d'une part le groupe fondamental d'une vari\'et\'e irr\'eductible orientable
est sans torsion par le th\'eor\`eme de la sph\`ere 
(voir ci-dessus la preuve du th\'eor\`eme 12),
et d'autre part le groupe fondamental d'une vari\'et\'e r\'eductible orientable
est un produit libre d'un groupe sans torsion et de groupes finis
(produit libre provenant d'une d\'ecomposition de Kneser-Milnor).
Par suite, si le groupe fondamental d'une $3$-vari\'et\'e orientable est infini,
il poss\`ede des \'el\'ements d'ordre infini~; 
on v\'erifie la m\^eme assertion pour le cas d'une $3$-vari\'et\'e non orientable
en consid\'erant le sous-groupe correspondant au rev\^etement d'orientation. 

  En particulier, 
$\Gamma$ poss\`ede  un sous-groupe 
isomorphe \`a $\Bbb Z \times \Bbb Z /2\Bbb Z$.
Un r\'esultat d'Epstein d\'ej\`a cit\'e (th\'eor\`eme 9.12 de \cite{Hemp--76})
implique que $M$ est de la forme $M_1 \sharp\, R$, 
o\`u $R$ est une vari\'et\'e ferm\'ee non orientable
telle que $\pi_1(R) = \Bbb Z \times \Bbb Z /2\Bbb Z$.
De plus $M_1$ est simplement connexe (proposition 3.i),
de sorte que $\Gamma =  \Bbb Z \times \Bbb Z /2\Bbb Z$.
Le th\'eor\`eme 12.10 de \cite{Hemp--76} montre que
$\hat M$ a le type d'homotopie de $\Bbb P^2 \times \Bbb S^1$.

\smallskip

R\'eciproquement, si $\hat M$ a le type d'homotopie de 
$\Bbb P^2 \times \Bbb S^1$, alors 
$\Gamma = \Bbb Z \times \Bbb Z /2\Bbb Z$ n'est pas cci~;
de m\^eme, le groupe di\'edral infini n'est pas cci~;
enfin, si $\Cal P(\tilde M)$ est une vari\'et\'e de Seifert,
$\Gamma$ n'est pas cci par la proposition 2 et le lemme 4.
\hfill$\square$
\enddemo

\medskip

\proclaim{17.\ Proposition}
Soient $M$ une $3$-vari\'et\'e non orientable $\Bbb P^2$-irr\'eductible
et $\Gamma$ son groupe fondamental.
Alors $\Gamma$ n'est pas cci si et seulement si $M$ est une vari\'et\'e de Seifert.
\endproclaim

\demo{D\'emonstration} 
Si $M$ est de Seifert, $\Gamma$ n'est pas cci par la proposition 2.

Supposons que $\Gamma$ n'est pas cci.
Le groupe $\Gamma$ est infini sans torsion par un r\'esultat d'Epstein d\'ej\`a
invoqu\'e (th\'eor\`eme 9.8 de \cite{Hemp--76}).
Le lemme 16 implique donc que $\Cal P (\tilde M)$ est une vari\'et\'e de Seifert~;
de plus $\tilde M$ est irr\'eductible (lemme 10.4 de \cite{Hemp--76}),
de sorte que $\tilde M = \Cal P (\tilde M)$ est une vari\'et\'e de Seifert.
Ceci ach\`eve la preuve gr\^ace au r\'esultat suivant.

{\it
Une vari\'et\'e non orientable irr\'eductible dont le rev\^etement d'orientation
est une vari\'et\'e de Seifert est elle-m\^eme une vari\'et\'e de Seifert.
}
C'est une cons\'equence dans le cas sans bord
du th\'eor\`eme~3 de \cite{Whit--92}
et dans le cas avec bord 
du r\'esultat principal de \cite{Toll--78}.
\phantom{aaaa}\hfill $\square$
\enddemo

\bigskip
\head{\bf 
X.~Le th\'eor\`eme principal 
pour les vari\'et\'es non orientables 
}\endhead
\medskip

\proclaim{18.\ Proposition} 
Soit $M$ une $3$-vari\'et\'e non orientable irr\'eductible
qui ne contient pas de $\Bbb P^2 \times I$ non standard
et dont le groupe fondamental $\Gamma$ est infini.
Alors $\Gamma$ n'est pas cci 
si et seulement si $M$ est une vari\'et\'e de Seifert modulo $\Bbb P$.
\endproclaim

\demo{D\'emonstration}
Si $M$ est de Seifert  modulo $\Bbb P$, 
le groupe $\Gamma$ n'est pas cci par
le th\'eor\`eme~15.

\medskip 

Supposons que $\Gamma$ n'est pas cci.
Vu la proposition 17, nous pouvons supposer que $M$
contient une copie $P$ du plan projectif, plong\'e, bilat\`ere~;
son groupe fondamental s'injecte dans celui de $M$
(remarque 14.i).
Pour la suite de la preuve, nous distinguons plusieurs cas.

\medskip

  {\it Cas (A)~: le plan projectif $P$ n'est pas parall\`ele au bord de $M$.}

   {\it Cas (A.1)~: de plus, $P$  s\'epare $M$ en deux vari\'et\'es $M_1,M_2$.}
Alors $M_1$ et $M_2$ sont non orientables 
(car \`a bord contenant un plan projectif)
et irr\'eductibles  (car $M$ est irr\'eductible).
De plus les groupes fondamentaux 
$\Gamma_1 = \pi_1(M_1),\Gamma_2 = \pi_1(M_2)$ sont infinis~;
en effet,  si $\Gamma_1$ \'etait fini, 
$M_1$ aurait le type d'homotopie d'un $\Bbb P^2 \times I$ 
par le th\'eor\`eme 9.6 de \cite{Hemp--76},
donc $M_1$ serait hom\'eomorphe \`a $\Bbb P^2 \times I$
vu les hypoth\`eses,
et donc $P$ serait parall\`ele au bord.

Le groupe $\Gamma$ est un produit amalgam\'e $\Gamma_1 \ast_{C_2} \Gamma_2$.
Par la proposition 3.iii, 
le groupe $C_2$ est normal dans $\Gamma_1$ et dans $\Gamma_2$.
Les groupes infinis $\Gamma_1$ et $\Gamma_2$ contiennent donc chacun
un sous-groupe isomorphe \`a $\Bbb Z \times C_2$.
Le corollaire 4.2 de \cite{Swar--73} 
implique que $M_1$ et $M_2$ sont homotopes \`a $\Bbb P^2 \times \Bbb S^1$,
ce qui est absurde car $M_1$ et $M_2$ sont \`a bord.

Le cas (A.1) n'est donc pas possible.

   {\it Cas (A.2)~: le plan $P$ ne s\'epare pas la vari\'et\'e $M$.}
Le groupe $\Gamma$ est donc une extension HNN
de base un groupe $\Gamma_1$ 
au-dessus d'un groupe $C_2$ d'ordre deux.
Si $C_2 = \Gamma_1$, alors
$\Gamma \approx \Bbb Z \times C_2$ 
et $M$ a le type d'homotopie de $\Bbb P^2 \times \Bbb S^1$,
donc $M$ est hom\'eomorphe \`a $\Bbb P^2 \times \Bbb S^1$
qui est une vari\'et\'e de Seifert.

On peut donc supposer que $C_2$ est un sous-groupe propre de $\Gamma_1$.
Alors $C_2$ est normal dans $\Gamma_1$
(sinon $\Gamma$ serait cci par la proposition 3.v).
Donc $\Gamma$ poss\`ede un sous-groupe isomorphe \`a $\Bbb Z \times C_2$.
Le r\'esultat de Swarup d\'eja cit\'e implique que $M$ 
est hom\'eomorphe \`a $\Bbb P^2 \times \Bbb S^1$,
qui est une vari\'et\'e de Seifert.

\medskip

{\it Cas (B)~: tout plan projectif plong\'e dans $M$ est parall\`ele au bord.}
Notons $\sigma$ l'involution du rev\^etement d'orientation $\tilde M$ de $M$
et $p : \tilde M \longrightarrow M$ la projection de rev\^etement.

\medskip

{\it Affirmation (I)~:} $\tilde M$ ne contient pas de boule non standard.

Supposons (ab absurdo) que $\tilde M$ contienne une boule non standard, 
dont nous notons $S$ le bord. 
On peut supposer que l'une des deux situations suivantes est r\'ealis\'ee~:
$\sigma(S) \cap S = \emptyset$, ou $\sigma(S) = S$
(voir la preuve du lemme 1 de \cite{Toll--70},
ou la preuve du lemme 2 de \cite{HeWhi--94}).

{\it Cas o\`u $\sigma(S) \cap S = \emptyset$.} 
Alors $M$ contient aussi une boule non standard, 
ce qui est absurde car $M$ est irr\'eductible 
et non simplement connexe (parce que non orientable).

{\it Cas o\`u $\sigma(S) = S$.}
Alors $p(S)$ est un plan projectif plong\'e dans $M$,
donc $p(S)$ est parall\`ele au bord dans $M$.
Il en r\'esulte que $S$ est parall\`ele au bord dans $\tilde M$,
et par suite que $\tilde M$ est obtenue par recollement d'une boule non standard
avec un produit $\Bbb S^2 \times I$~;
en particulier $\tilde M$ est simplement connexe, 
ce qui est absurde car $\Gamma = \pi_1(M)$ est infini.

L'affirmation (I) est ainsi d\'emontr\'ee,
de sorte que $\Cal P (\tilde M) = \overline{M}$,
o\`u nous \'ecrivons $\overline{M}$ pour $\hat{\tilde M}$
(la signification du chapeau est celle d\'efinie au \S~VI).
Le lemme 16 montre que ou bien $\pi_1(\tilde M)$ est di\'edral infini,
ou bien  $\overline{M}$ est une vari\'et\'e de Seifert.

\medskip

{\it Affirmation (II)~: $\pi_1(M)$ contient un sous-groupe cyclique distingu\'e
non r\'eduit \`a un \'el\'ement.}

{\it Cas o\`u $\pi_1(\tilde M)$ est un groupe di\'edral infini.}
Le sous-groupe sans torsion d'indice deux du groupe di\'edral infini 
est caract\'eristique, 
donc $\pi_1(M)$ contient un sous-groupe cyclique infini distingu\'e.

{\it Cas o\`u $\overline{M}$ est une vari\'et\'e de Seifert.} 
Alors, le groupe fondamental de $\overline{M}$
contient un sous-groupe normal infini, 
engendr\'e par les fibres r\'eguli\`eres. 
La vari\'et\'e de Seifert orientable $\overline{M}$
est ou bien hom\'eomorphe \`a 
l'une de $\Bbb S^2 \times \Bbb S^1$, $\Bbb P^3 \sharp\, \Bbb P^3$,
ou bien irr\'eductible
(lemme VI.7 de \cite{Jaco--77}). Dans le premier cas,
$\pi_1(M)$ contient \'evidemment un sous-groupe cyclique infini normal.
On peut donc supposer que $\overline{M}$ est de plus irr\'eductible.

L'involution d'orientation $\sigma$ sur $\tilde{M}$
s'\'etend en une involution $\overline{\sigma}$ sur $\overline{M}$
qui a un nombre fini de points fixes 
(un par composante sph\'erique de $\partial\tilde{M}$).
Notons $X$ l'espace des orbites de la fibration de Seifert sur $\overline{M}$.
L'orbi\'et\'e $X$ ne peut pas \^etre $\Bbb S^2$ 
avec $0$, $1$ ou $2$ points coniques 
(sinon $\overline{M}$ serait $\Bbb S^2 \times \Bbb S^1$ ou un espace lenticulaire,
ce qui est exclu car $\overline M$ est irr\'eductible 
et $\pi_1(\overline M)$ est infini),
ni $\Bbb S^2$ avec $3$ points coniques
(par un argument de la preuve du th\'eor\`eme~2 dans \cite{HeWh--94}).
Nous pouvons donc appliquer le th\'eor\`eme principal de \cite{Toll--78},
qui nous assure que $\overline{\sigma}$ pr\'eserve 
la fibration de Seifert de $\overline{M}$.

Notons $t$ un g\'en\'erateur du sous-groupe cyclique infini normal
de $\pi_1(\tilde M)$, qui s'identifie \`a $\pi_1(\overline{M})$,
et $a$ un \'el\'ement de $\pi_1(M)$ qui n'est pas dans 
le sous-groupe d'indice deux $\pi_1(\tilde M)$.
Nous affirmons que $tat^{-1}$ est $t$ ou $t^{-1}$,
de sorte que le sous-groupe engendr\'e par $t$ 
est encore normal dans $\pi_1(M)$.

En effet, choisissons un point base $x_0 \in M$ 
et un relev\'e $\tilde x_0 \in \tilde M$
de telle sorte que $\tilde x_0$ 
soit dans une fibre r\'eguli\`ere de $\overline{M}$.
L'inclusion $\tilde M \subset \overline{M}$
fournit une identification 
$\pi_1(\tilde M, \tilde x_0) = \pi_1(\overline{M}, \tilde x_0)$
et la projection de rev\^etement $p$ une identification de 
$\pi_1(\tilde M, \tilde x_0)$ 
avec un sous-groupe d'indice deux dans $\pi_1(M, x_0)$.
Choisissons encore 
un lacet $\alpha$ dans $M$ bas\'e en $x_0$ 
repr\'esentant $a$,
dont le relev\'e $\tilde \alpha$ dans $\tilde M$
connecte $\tilde x_0$ \`a $\sigma(\tilde x_0)$,
un lacet $\tilde \tau$ dans $\tilde M$ bas\'e en $\tilde x_0$
qui repr\'esente $t$,
et notons $\overline{\tau}$ la fibre de $\overline{M}$ contenant $\tilde x_0$
orient\'ee de telle sorte qu'elle repr\'esente aussi $t$
(en particulier, $\tilde \tau$ et $\overline{\tau}$ 
sont homotopes dans $\overline{M}$).

Consid\'erons le lacet 
$\tilde \alpha \sigma(\tilde \tau) (\tilde \alpha)^{-1}$
de $\tilde M$.
D'une part, il se projette dans $M$ 
sur un lacet $\alpha \tau \alpha^{-1}$
qui repr\'esente $ata^{-1}$.
D'autre part, il est librement homotope dans $\overline{M}$
au lacet $\overline{\sigma}(\overline{\tau})$,
qui est une fibre r\'eguli\`ere orient\'ee, 
donc qui est aussi librement homotope dans $\overline{M}$
\`a $\overline{\tau}$ ou \`a $(\overline{\tau})^{-1}$~;
par suite $ata^{-1}$ est conjugu\'e dans $\pi_1(\overline{M},\tilde x_0)$
\`a $t$ ou \`a $t^{-1}$~;
mais $t$ engendre dans $\pi_1(\overline{M},\tilde x_0)$ 
un sous-groupe qui est cyclique infini normal, 
donc $ata^{-1}$ est l'un de $t$ ou $t^{-1}$,
comme affirm\'e plus haut.

Ceci ach\`eve la preuve de l'affirmation (II). Le th\'eor\`eme 15 permet de
conclure.
\hfill $\square$
\enddemo

\medskip

\proclaim{19.\ Th\'eor\`eme}
Soit $M$ une $3$-vari\'et\'e non orientable
qui ne contient pas de $\Bbb P^2 \times I$ non standard,
dont le groupe fondamental $\Gamma$ est infini et n'est pas di\'edral infini.
Alors $\Gamma$ n'est pas cci  si et seulement si 
$\Cal P (M)$ est une vari\'et\'e de Seifert modulo $\Bbb P$.

En particulier, pour un groupe fondamental infini $\Gamma$
d'une $3$-vari\'et\'e non orientable,
les quatre propri\'et\'es suivantes sont \'equivalentes~:
\roster
\item"(i)" $\Gamma$ est groupe fondamental 
   d'une vari\'et\'e de Seifert modulo $\Bbb P$~; 
\item"(ii)" $\Gamma$ poss\`ede un sous-groupe normal cyclique infini~;
\item"(iii)" $\Gamma$ n'est pas cci~;
\item"(iv)" $\Gamma$ poss\`ede un sous-groupe d'indice $2$
qui est groupe fondamental d'une vari\'et\'e de Seifert orientable.
\endroster
Lorsque $\Gamma$ n'a pas d'\'el\'ement d'ordre deux,
ces propri\'et\'es sont encore \'equivalentes \`a~:
\roster
\item"(v)" $\Gamma$ est groupe fondamental 
   d'une vari\'et\'e de Seifert.
\endroster
\endproclaim

\demo{D\'emonstration}
Le th\'eor\`eme 19 r\'esulte de la proposition 18
comme le th\'eor\`eme 13 r\'esulte du th\'eor\`eme  12.
\hfill $\square$
\enddemo

\bigskip
\head{\bf 
XI.~Cas d'un entrelacs et cas d'une vari\'et\'e hyperbolique 
}\endhead
\medskip

\proclaim{20.\ Corollaire} 
Soit $L$ un noeud dans $\Bbb S^3$. Alors
le groupe de $L$ est cci si et seulement si $L$ n'est pas un noeud
torique.
\endproclaim

\demo{D\'emonstration} D'une part, le groupe d'un noeud  est
infini (le th\'eor\`eme de Hurwitz et la dualit\'e d'Alexander
montrent que son ab\'elianis\'e est cyclique infini). 
D'autre part, on sait essentiellement depuis la classification par Seifert 
des structures de Seifert sur la sph\`ere que
le compl\'ementaire d'un noeud dans $\Bbb S^3$ 
poss\`ede une fibration de Seifert 
si et seulement si le noeud est torique 
(voir par exemple \cite{Mose--71}).
\hfill $\square$
\enddemo

\demo{Remarque} Plus g\'en\'eralement, soit $L$ un entrelacs dans $\Bbb S^3$. 
Les conditions suivantes sont \'equivalentes \cite{BuMu--70}~:
\roster
\item"(i)" le centre du groupe fondamental $\pi_1(\Bbb S^3 \setminus L)$ n'est pas
      r\'eduit \`a $1$~;
\item"(ii)" le groupe fondamental $\pi_1(\Bbb S^3 \setminus L)$ n'est pas cci~;
\item"(iii)" l'entrelacs $L$ est une r\'eunion finie de fibres d'une fibration de
Seifert de $\Bbb S^3$.
\endroster
\enddemo

\medskip

Voici enfin un \'enonc\'e 
qu'on pourrait sans doute d\'eduire  de ce qui pr\'ec\`ede, 
mais dont nous pr\'ef\'erons donner une d\'emonstration 
dans un autre contexte.

\medskip

\proclaim{21.\ Proposition} Le groupe fondamental d'une
vari\'et\'e hyperbolique orientable de volume fini 
est toujours un groupe cci.
\endproclaim

\demo{D\'emonstration} 
Un tel groupe fondamental est un r\'eseau
dans la composante connexe du groupe des isom\'etries de l'espace hyperbolique.
Or c'est une cons\'equence du
th\'eor\`eme de densit\'e de Borel \cite{Bore--60} que, 
plus g\'en\'eralement, 
dans un groupe de Lie $G$ connexe simple non compact 
\`a centre r\'eduit \`a un \'el\'ement,
tout r\'eseau $\Gamma$ est cci.

Plus pr\'ecis\'ement, soit $\gamma \in \Gamma$, $\gamma \ne 1$. 
Il s'agit de montrer qu'il existe une suite infinie $(\gamma _j)_{j \ge 1}$ 
telle que les conjugu\'es $\gamma_j \gamma \gamma_j^{-1}$ 
soient distincts deux \`a deux.
On pose $\gamma _1 = 1$ et on proc\`ede par r\'ecurrence. 

Supposons qu'il existe une suite $(\gamma _j)_{1 \le j \le k}$ telle que
les \'el\'ements $\gamma _j^{-1} \gamma \gamma _j$ ($1 \le j \le k)$
sont distincts deux \`a deux. 
Pour chaque $j \in \{1,...,k\}$, le ferm\'e de Zariski
$$
    Z_G(\gamma) \gamma _j \, = \,
    \{ g \in G : g^{-1} \gamma g = \gamma _j^{-1}\gamma \gamma _j \}
$$
est distinct de $G$, sinon $\gamma$ serait dans le centre de $G$.
(La notation $Z_G(\gamma)\gamma _j$ indique le translat\'e par
$\gamma _j$ d'un centralisateur dans $G$.) 
Comme la r\'eunion d'un nombre fini de ferm\'es
de Zariski distincts de $G$ est distincte de $G$, 
la propri\'et\'e de $\Gamma$ d'\^etre Zariski-dense dans $G$
implique qu'il existe $\gamma_{k+1} \in \Gamma$ tel que 
$$
    \gamma_{k+1} \, \notin \, 
   \bigcup _{1 \le j \le k}  Z_G(\gamma)\gamma _j .
$$
Par suite les \'el\'ements $\gamma _j^{-1} \gamma \gamma _j$ 
($1 \le j \le k+1)$ sont distincts deux \`a deux.

Il en r\'esulte que la classe de conjugaison de $\gamma$ dans $\Gamma$ est
infinie.
\hfill $\square$ 
\enddemo

Notons que la conclusion de la proposition et sa d\'emonstration  valent 
pour un sous-groupe Zariski-dense 
d'un groupe de Lie alg\'ebrique r\'eel de centre r\'eduit \`a $1$,
et en parti\-culier pour une vari\'et\'e hyperbolique 
de volume fini en toute  dimension $n \ge 2$. 

\medskip

Plus g\'en\'eralement, on peut montrer que le groupe fondamental d'une
$3$-vari\'et\'e (compacte ou non) hyperbolique non \'el\'ementaire 
est un groupe cci.
Une vari\'et\'e hyperbolique, quotient de l'espace hyperbolique $\Bbb H^3$
par un groupe discret sans torsion $\Gamma$ d'isom\'etries de $\Bbb H^3$,
est dite {\it \'el\'ementaire} si l'ensemble limite de $\Gamma$ dans le bord
$\Bbb S^2$ de $\Bbb H^3$ a au plus deux points, c'est-\`a-dire si $\Gamma$
poss\`ede un sous-groupe d'indice fini isomorphe \`a $\Bbb Z^d$ ($d \le 2$).

\bigskip
\head{\bf 
XII.~Groupes de dimension cohomologique trois \`a dualit\'e de Poincar\'e
}\endhead
\medskip

Soit $\Gamma$ un groupe 
de dimension cohomologique trois \`a dualit\'e de Poincar\'e~;
nous \'ecrirons succintement~: \lq\lq soit $\Gamma$ {\it un groupe DP(3)}\rq\rq ,
et nous renvoyons au \S~VIII.10 de \cite{Brow--82} pour la d\'efinition.
Rappelons toutefois 
qu'un groupe DP(3) est sans torsion,
que tout groupe fondamental 
d'une $3$-vari\'et\'e ferm\'ee asph\'e\-rique
\footnote{
Ou, de mani\`ere \'equivalente, tout groupe fondamental
d'une $3$-vari\'et\'e ferm\'ee $\Bbb P^2$-irr\'eductible
\`a rev\^e\-tement universel non compact.
}
est un groupe DP(3), 
et qu'on ne conna\^{\i}t pas d'autre exemple.
Notons aussi que, si le groupe fondamental d'une $3$-vari\'et\'e $M$ est DP(3),
alors la vari\'et\'e $\Cal P(M)$ est n\'ecessairement 
ferm\'ee et $\Bbb P^2$-irr\'eductible

L'analogue de la conjecture des fibr\'es de Seifert pour les groupes DP(3)
est le corollaire 0.5 de \cite{Bowd--04}~:

\medskip

\proclaim{22.\ Proposition (Bowditch)}
Un groupe DP(3) qui contient un sous-groupe normal cyclique infini
est le groupe fondamental d'une vari\'et\'e de Seifert ferm\'ee.
\endproclaim

\medskip

Dans le cadre de ce paragraphe, voici l'analogue des th\'eor\`emes 13 et 19.

\medskip

\proclaim{23.\ Proposition}
Un groupe DP(3) n'est pas cci si et seulement si
c'est le groupe fondamental d'une vari\'et\'e de Seifert.
\endproclaim

\demo{D\'emonstration} Soit $\Gamma$ un groupe DP(3) qui n'est pas cci~;
il s'agit de montrer qu'il existe une vari\'et\'e de Seifert 
dont $\Gamma$ est le groupe fondamental.
Par le lemme 9, il existe dans $\Gamma$ un sous-groupe $\Gamma_0$ 
qui est distingu\'e, ab\'elien, de type fini et infini~;
de plus, sa dimension cohomologique est au plus trois.
Par suite, $\Gamma_0$ est isomorphe \`a l'un des groupes
$\Bbb Z$, $\Bbb Z^2$, $\Bbb Z^3$.

Supposons d'abord $\Gamma_0$ d'indice fini.
Sa dimension cohomologique est alors trois
(th\'eo\-r\`eme VIII.3.1 de \cite{Brow--82}),
de sorte que $\Gamma_0$ est isomorphe \`a $\Bbb Z^3$.
L'argument du quatri\`eme cas de la preuve du th\'eor\`eme~12
permet de conclure.

On suppose d\'esormais $\Gamma_0$ d'indice infini,
et donc isomorphe \`a l'un des groupes $\Bbb Z$, $\Bbb Z^2$,
par le th\'eor\`eme~1 de \cite{Hill--87}.
Si $\Gamma_0 = \Bbb Z$, la proposition 22 permet de conclure.

On peut donc supposer que $\Gamma_0$ est isomorphe \`a $\Bbb Z^2$.
Le corollaire de \cite{Hill--87} implique alors qu'il existe
une $3$-vari\'et\'e $\Bbb P^2$-irr\'eductible 
dont $\Gamma_0$ est le groupe fondamental.
Le th\'eor\`eme 12 et la proposition 17 permettent de conclure.
\hfill$\square$
\enddemo

\bigskip
\head{\bf
XIII.~La propri\'et\'e cci forte 
\\
Une application \`a certaines repr\'esentations unitaires
}
\endhead
\medskip

   Dans l'introduction, nous avons rappel\'e que
la propri\'et\'e pour un groupe $\Gamma$ d'\^etre cci 
se traduit en termes de l'alg\`ebre de von Neumann $W^*_{\lambda}(\Gamma)$~;
elle a aussi des cons\'equences 
pour l'\'etude des repr\'esentations unitaires de $\Gamma$, 
et donc \cite{Harp} des C$^*$-alg\`ebres associ\'ees.
Avant d'en citer une, 
rappelons d'abord une d\'efinition de \cite{BeHa--94}~:
un groupe $\Gamma$ est {\it fortement cci} si, 
pour toute partie finie $F$ de $\Gamma$ disjointe de $\{1\}$, 
il existe une suite infinie $(\gamma _j)_{j \ge 1}$
d'\'el\'ements de $\Gamma$ telle que, pour tout $x \in F$, 
les \'el\'ements  $\gamma _j^{-1}x\gamma _j$ sont distincts deux \`a deux. 
En fait, cette propri\'et\'e est \'equivalente \`a la propri\'et\'e cci
(comme d\'ej\`a not\'e dans une note de \cite{BeHa--94} ajout\'ee aux \'epreuves).

\medskip

\proclaim{24.~Proposition} 
Un groupe est cci si et seulement s'il est fortement cci.
\endproclaim

\demo{D\'emonstration} Soit $\Gamma$ un groupe non r\'eduit \`a un \'el\'ement,
et cci~;  il s'agit de montrer que $\Gamma$ est fortement cci.
Soit $F \subset \Gamma \setminus \{1\}$ un ensemble fini.
Nous affirmons que, pour tout $k \ge 1$, 
il existe une suite $\gamma_1,\hdots,\gamma_k \in \Gamma$
telle que, pour tout $f \in F$, 
les \'el\'ements $\gamma_j f \gamma_j^{-1}$ ($j=1,\hdots,k$) 
sont distincts deux \`a deux.
On pose $\gamma_1=1$ et on proc\`ede par r\'ecurrence sur $k$,
en supposant l'affirmation d\'emontr\'ee jusqu'\`a $k$.

Comme $\Gamma$ est cci, le centralisateur $Z_{\Gamma}(f)$ de $f$ dans $\Gamma$ 
est d'indice infini pour tout $f \in F$.
Or c'est un fait classique qu'un groupe infini (ici $\Gamma$) 
n'est jamais r\'eunion d'un nombre fini de classes \`a droite
suivant des sous-groupes d'indices infinis (ici les $Z_{\Gamma}(f) \gamma_j$)~;
voir le lemme 4.1 de \cite{Neum--54}.
Nous pouvons donc choisir $\gamma_{k+1} \in \Gamma$ tel que
$\gamma_{k+1} \notin \cup Z_{\Gamma}(f) \gamma_j$
(r\'eunion sur $f \in F$ et $j = 1,\hdots,k$).
Il est alors imm\'ediat de v\'erifier que,
$\gamma_j f \gamma_j^{-1} \ne \gamma_{k+1}f\gamma_{k+1}^{-1}$
pour tous $j \in \{1,\hdots,k\}$ et $f \in F$.
\hfill $\square$
\enddemo 

Rappelons encore qu'une repr\'esentation unitaire $\pi$ 
d'un groupe $\Gamma$ dans un espace de Hilbert $\Cal H_{\pi}$ 
est dite {\it de classe $(\Cal C _0)$} 
si, pour tout vecteur unit\'e $\xi \in \Cal H _{\pi}$,
la fonction de type positif $\gamma \mapsto <\xi | \pi (\gamma)\xi >$ 
tend vers $0$ \`a l'infini de $\Gamma$.
Nous \'ecrivons $\pi \prec \rho$ si une repr\'esentation unitaire $\pi$
est {\it faiblement contenue} dans une repr\'esentations unitaire~$\rho$.

\medskip

\proclaim{25.~Proposition}
Si $\Gamma$ est un groupe cci 
et si $\pi$ une repr\'esentation de $\Gamma$ de classe $(\Cal C _0)$, 
alors $\lambda _{\Gamma} \prec \pi$.
\endproclaim

\demo{D\'emonstration} Soit $\delta : \Gamma \to \Bbb C$ 
la fonction d\'efinie par
$\delta (\gamma ) = <\xi _1 | \lambda _{\Gamma}(\gamma ) \xi _1 >$ 
o\`u  $\xi _1 \in l^2(\Gamma )$ est la fonction caract\'eristique de $\{1\}$; 
nous avons donc $\delta (1) = 1$ et $\delta (\gamma ) = 0$ si $\gamma \ne 1$.
Vu la proposition 18.1.4 de \cite{DiC$^*$},
il suffit de v\'erifier que $\delta$ est approchable 
par des fonctions de type positif associ\'ees \`a $\pi$. \par

   Soit $F$ une partie finie de $\Gamma$ disjointe de $\{1\}$. 
Soit $(\gamma _j)_{j \ge 1}$ une suite 
comme dans la d\'efinition de \lq\lq fortement cci\rq\rq . 
Par hypoth\`ese sur $\pi$, nous pouvons
choisir un vecteur unit\'e $\xi _0 \in \Cal H _{\pi}$ tel que
$$
   \lim _{j \to \infty} <\xi _0 | \pi (\gamma _j^{-1}x\gamma _j )\xi _0 > = 
   0  \qquad \text{pour tout} \quad f \in F
$$
car $\gamma _j^{-1} f \gamma _j \ne \gamma _k^{-1} f \gamma _k$ si $j \ne k$. 
D\'efinissons $\phi _{F,j} : \Gamma \to 	\Bbb C$ par
$$
   \phi _{F,j}(\gamma) \, = \,   
   <\pi (\gamma _j)\xi _0 | \pi (\gamma ) \pi(\gamma _j) \xi_0 > .
$$
Nous avons donc bien $\lim_{j\to\infty}\phi_{F,j}(f) = \delta (f)$ pour tout
$f \in F \cup \{1\}$. 
\hfill $\square$ \enddemo

\medskip

En d'autres termes, et avec les notations de \cite{Harp},
si $\Gamma$ est un groupe cci 
et si $\pi$ est une repr\'esentation de $\Gamma$ de classe $(\Cal C _0)$, 
alors la C$^*$-alg\`ebre r\'eduite $C^*_{\lambda}(\Gamma)$
est naturellement un quotient de la C$^*$-alg\`ebre $C^*_{\pi}(\Gamma)$
associ\'ee \`a $\pi$.

\medskip

C'est une question naturelle que de demander dans quel cas,
pour un $3$-groupe $\Gamma$ qui est cci,
la C$^*$-alg\`ebre r\'eduite $C^*_{\lambda}(\Gamma)$ est simple.
(C'est par exemple toujours le cas si $\Gamma$
est le groupe fondamental d'une vari\'et\'e hyperbolique close.)

\enddocument

\head
Vieille chute 1 (\`a garder ?)
\endhead

Le groupe fondamental d'une vari\'et\'e hyperbolique orientable
de volume fini est cci
\footnote{
Dire pourquoi en moins de 3 lignes ?????
\newline
Remarque : une vari\'et\'e dont le rev\^etement universel est contractile est
toujours irr\'eductible (r\'esulte du th\'eor\`eme de la sph\`ere).
\newline
Par ailleurs, on connait des conditions suffisantes garantissant 
l'existence d'une structure hyperbolique de volume fini
sur une $3$-vari\'et\'e sans bord 
(?????????????)
ou sur l'int\'erieur d'une $3$-vari\'et\'e \`a bord non vide
(th\'eor\`eme B de \cite{Morg--84}).
\newline
??? Voir peut-\^etre aussi le th\'eor\`eme 2.3 dans Thurston,
Bull. AMS 6 (1982) 357--381.
}
\hskip-.1cm
.

\head
Vieille chute 2 (\`a garder ?)
\endhead

\noindent !!!!!!! La suite \`a v\'erifier !!!!!!!

\noindent
Une vari\'et\'e de Seifert $M$ qui n'est pas Haken peut \^etre
\roster
\item
un espace lenticulaire,
\item
$\Bbb P^3 \sharp\, \Bbb P^3$,
\item
un fibr\'e de Seifert au-dessus de $\Bbb S^2$ avec exactement $3$
fibres exceptionnelles, tel que $H_1(M,\Bbb Z)$ est infini. (???
Corriger~: \lq\lq fini\rq\rq \ ???)
\endroster
(Th\`ese de JPP, page 120.)

\head
Stallings et Epstein
\endhead

(ii) La tradition orale (que nous avons entendue de David Epstein)
semble indiquer que cette conjecture ait \'et\'e formul\'ee par Waldhausen
\`a la suite d'un r\'esultat de Stallings selon lequel,
sous des hypoth\`eses convenables, une $3$-vari\'et\'e
dont le groupe fondamental est un produit de deux groupes infinis
est n\'ecessairement produit d'une surface et d'un cercle.
\newline
????????????? A repenser !!!!!!!!!!!!!

\head
Anciennement apr\`es le lemme 15
\endhead

\medskip

\proclaim{***.~Corollaire}
Soit $\Gamma$ un groupe infini
qui est groupe fondamental d'une $3$-vari\'et\'e non orientable. 
Les trois conditions suivantes sont \'equivalentes~:
\roster
\item"(i)" $\Gamma$ poss\`ede un sous-groupe d'indice $2$ 
   qui est groupe fondamental d'une vari\'et\'e de Seifert~; 
\item"(ii)" $\Gamma$ poss\`ede  un sous-groupe d'indice $2$ 
   qui poss\`ede un sous-groupe normal cyclique infini~;
\item"(iii)" $\Gamma$ n'est pas cci.
\endroster
\endproclaim

\demo{D\'emonstration} 
Le th\'eor\`eme 15 montre que (iii) implique (i).
Le lemme 1 montre que (i) implique (ii).
Vu qu'un groupe contenant un sous-groupe normal cyclique infini
n'est pas cci, le lemme 4 montre que (ii) implique (iii).
\hfill$\square$
\enddemo

\head 
Z.\ Pense-b\^ete pour PH
\endhead

\medskip
\subhead 
Id\'ee de la preuve pour 
$\quad \Gamma \  \text{cci}
\quad \Longleftrightarrow \quad 
W^*_{\lambda}(\Gamma) \ \text{facteur}$
\endsubhead
Rappel~: \`a cette \'epoque, le groupe \`a un \'el\'ement \'etait consid\'er\'e
comme \'etant cci !!!
\medskip

Soient $\varphi \in \Bbb C [\Gamma]$, $\gamma \in \Gamma$, 
$\xi \in \ell^2(\Gamma)$ et $x \in \Gamma$. 
Alors
$$
\aligned 
\left( 
\lambda_{\Gamma} (\gamma^{-1})
\lambda_{\Gamma}(\varphi) 
\lambda_{\Gamma}(\gamma) \xi
\right) (x) 
\, &= \, 
\left( 
\lambda_{\Gamma}(\varphi) 
\lambda_{\Gamma}(\gamma) \xi 
\right) (\gamma x) 
\, = \, 
\sum_{y,z \in \Gamma, yz = \gamma x} 
\lambda_{\Gamma} (\varphi) (y) \xi (\gamma^{-1} z)
\\
\, &= \, 
\sum_{y,v \in \Gamma, y\gamma v = \gamma x}
\lambda_{\Gamma}(\varphi) (y) \xi (v) .
\endaligned
$$
En particulier, si $\xi = \delta_e$~:
$$
\left(
\lambda_{\Gamma} (\gamma^{-1}) 
\lambda_{\Gamma}(\varphi)
\lambda_{\Gamma}(\gamma) \delta_e
\right) (x) 
\, = \, 
\left( \lambda_{\Gamma}(\varphi) \right) 
(\gamma x \gamma^{-1}) .
$$
Il en r\'esulte que l'\'el\'ement $\lambda_{\Gamma}(\varphi)$ est
central si et seulement si la fonction $\varphi$ est constante sur
les classes de conjugaison de $\Gamma$.

Modulo quelques pr\'ecautions 
(voir le lemme 4.2.12 de \cite{Saka--71}, d\'ej\`a cit\'e),
ceci implique que le centre de $W^*_{\lambda}(\Gamma)$ 
est r\'eduit \`a $\Bbb C$ si et seulement si $\Gamma$ est cci.

\medskip
\subhead Divers
\endsubhead
\medskip

Soit $M$ une $3$-vari\'et\'e telle que $\pi_1(M)$ est de type
fini. Alors il existe une sous-vari\'et\'e compacte $M_c$ de $M$
telle que l'inclusion induise un isomorphisme $\pi_1(M_c) \approx
\pi_1(M)$~; voir \cite{Hempel, page 73}.

\medskip
\subhead Vari\'et\'es de Seifert
\endsubhead
\medskip

\par\noindent
Lemme VII.7 de Jaco : $\Bbb S^2 \times \Bbb S^1$, ou $\Bbb P^3
\sharp\, \Bbb P^3$, ou irr\'eductible.
\par\noindent
Exo 6.8 de Jaco : si $M$ Seifert irr\'eductible, alors $M$
$\partial$-irr\'eductible ssi $M \ne$ tore solide.
\par\noindent
Jaco VI.29 : $M$ Seifert et rev\^etement de $M$ Seifert. Mieux :
JS \S \ II.6 pages 36 et 38.
\par\noindent
Jaco VI.15 : toute $M$ Seifert ou bien irr\'eductible ou bien
courte liste.
\par
Peut-\^etre mieux, Jaco VI.7~: la liste courte est $S^1 \times
S^2$, $\Bbb P^3 \sharp\, \Bbb P^3$.
\par\noindent
Jaco-Shalen p. 23 : d\'ef. du sous-groupe canonique.
\par\noindent
Jaco-Shalen p. 25, si $\Gamma$ infini, son sous-groupe caract.
d'indice $\le 2$~; cf aussi II.4.5 p. 26.

Soit $M$ une vari\'et\'e de Seifert orientable \`a $\Gamma$
infini, et $h \in \Gamma$ l'\lq\lq \'el\'ement fibre\rq\rq . Lemme
II.4.2 de \cite{JaSh--79}~: $h$ est d'ordre infini. Voir aussi VI.24
dans \cite{Jaco} (caract\'erisation des Seifert parmi les Haken).

\medskip

Structure de Seifert sur $\Bbb P^3 \sharp\, \Bbb P^3$~: voir ca
comme un quotient de $\Bbb S^2 \times [0,1]$ (avec $1=0$) modulo
l'involution d\'efinie par $(z,t) \longmapsto (-z,t+\frac{1}{2})$.

\medskip

\proclaim{Th \cite{Sco--83?  }} Soit $M$ une $3$-vari\'et\'e close
orientable irr\'eductible \`a groupe fondamental infini. Si $M$
poss\`ede un rev\^etement fini qui est de Seifert, alors $M$ est
de Seifert.
\endproclaim

\medskip
\subhead 
Vari\'et\'es ind\'ecomposables et irr\'eductibles
\endsubhead
\medskip

Dans cette th\'eorie, on convient que $\Bbb S^3$ est une
vari\'et\'es ind\'ecomposable.

\bigskip

Ne pas confondre $M \sharp\, \Bbb S^3 = M$ et $M \sharp\, \Bbb B^3 = M
\setminus \text{(open $3$-ball)}$.

\bigskip

Sur les ordres $M \longmapsto \hat M \longmapsto \Cal P (M)$ et $M
\longmapsto \Cal P (M) \longmapsto \hat M$.

Si on passe d'abord par $\hat M$, il n'y a pas de $3$-cellules
dans la d\'ecomposition de Kneser-Milnor. Si on d\'ecompose
d'abord $M$ (au lieu de $\hat M$), alors il faut peut-\^etre
distinguer ceux des facteurs $M_j$ qui sont des $3$-cellules des
autres.

\medskip

\proclaim{5.\ Lemme (r\'edaction pr\'ec\'edente)} 
Si $H$ est un sous-groupe d'indice 2 de
$G$, et si $H$ est cci, alors soit $G$ est aussi cci, 
soit $G =H \times \Bbb Z / 2\Bbb Z$.
\endproclaim

\demo {D\'emonstration} Notons $\bar{G}$ la r\'eunion des classes
de conjugaison finies de $G$. Alors $\bar{G}$ est un sous-groupe
de $G$ (clairement si $u\in\bar{G}$ alors $u^{-1}\in\bar{G}$~;
pour montrer que la loi est interne il suffit de remarquer qu'un
conjugu\'e de $ab$ est le produit d'un conjugu\'e de $a$ et d'un
conjugu\'e de $b$), qui est clairement distingu\'e dans $G$.

Remarquons maintenant que si $H$ est un sous-groupe d'indice $k$
de $G$, et $K$ est un sous-groupe de $G$ d'ordre $>k$, alors
$H\cap K\not=\{ e\}$. En effet, puisque si deux \'el\'ements
$u,v\in K$ sont dans la m\^eme classe de $G/H$ alors $uv^{-1}\in
H$, n\'ecessairement avec le principe des tiroirs, il existe un
\'el\'ement non trivial de $K$ dans $H$.

Puisque $H$ est cci, $H\cap \bar{G}=\{ e\}$, et on d\'eduit de ce
qui pr\'ec\`ede que $\bar{G}$ est d'ordre au plus 2. Si
$\bar{G}=\{1,a\}$ est d'ordre 2, alors n\'ecessairement $a\in G-H$
est d'ordre 2 et central dans $G$. Ainsi,
$G=H\times\{1,a\}=H\times\Bbb Z / 2\Bbb Z$. Et si $\bar{G}=1$ est
d'ordre 1, alors $G$ est cci.
\hfill$\square$
\enddemo

\subhead Vieille version du lemme 1 et de la proposition 2\endsubhead

\proclaim{1.\ Lemme} Soit $\Gamma$ le groupe fondamental d'une
vari\'et\'e de Seifert $M$  orientable. Si $\Gamma$ est
infini, la classe d'homotopie d'une fibre r\'eguli\`ere de $M$
engendre dans $\Gamma$ un sous-groupe normal cyclique infini.
\endproclaim

\demo{D\'emonstration} Si $M$ est r\'eductible, $M$ est
hom\'eomorphe \`a $\Bbb S^2 \times \Bbb S^1$ ou \`a la somme
connexe $\Bbb P^3 \sharp\, \Bbb P^3$ de deux espaces projectifs
(lemme VI.7 de \cite{Jaco--77})~; vu que~ $\pi_1\left( \Bbb S^2
\times \Bbb S^1 \right) \approx \Bbb Z$ et $\pi_1\left( \Bbb P^3
\sharp\, \Bbb P^3 \right) \approx  D_{\infty}$ (le groupe di\'edral
infini), le lemme est \'evident dans ces cas.

Supposons donc $M$ irr\'eductible. 
Si $M$ est $\partial$-r\'eductible, 
alors $M$ est hom\'eomorphe \`a un tore
solide (exercice VI.8 de \cite{Jaco--77} ou th\'eor\`eme 10.4 de
\cite{FoMa--97}) et $\pi_1(M) \approx \Bbb Z$.

Supposons donc que $M$ est $\partial$-irr\'eductible. 
La conclusion du lemme ci-dessus fait
alors partie du lemme II.4.2 de \cite{JaSh--79}. 
\hfill $\square$
\enddemo

\medskip

\proclaim{2.\ Proposition} Le groupe fondamental d'une vari\'et\'e
de Seifert non simplement connexe poss\`ede une classe de
conjugaison finie autre que $1$.
\endproclaim

\demo{D\'emonstration} Si $\Gamma$ est comme au lemme 1, la classe
de conjugaison dans $\Gamma$ d'un \'el\'ement distinct de
l'identit\'e dans un sous-groupe normal cyclique infini contient
au plus deux \'el\'ement, et donc $\Gamma$ n'est pas cci. Il
suffit donc de consid\'erer le cas d'une vari\'et\'e de Seifert
$M$ non orientable \`a groupe fondamental $\Gamma$ infini. Notons
$\tilde \Gamma$ le sous-groupe d'indice deux de $\Gamma$
correpondant au rev\^etement d'orientation $p : \tilde M
\longrightarrow M$. La vari\'et\'e $\tilde M$ est aussi de Seifert
(lemme II.6.1 de \cite{JaSh--79}). Choisissons une fibre
r\'eguli\`ere $T$ de $M$, munie d'une orientation~; d\'esignons
par $t$ la classe de $T$ dans $\Gamma$ et par $\tilde t \in \tilde
\Gamma$ la classe dans $\tilde \Gamma$ d'une composante connexe de
$p^{-1}(T)$. Ainsi $\tilde t$ est \'egal soit \`a $t$ soit \`a
$t^2$. Dans tous les cas, la classe de conjugaison dans $\Gamma$
d'une puissance non nulle de $\tilde t$ contient au plus deux
\'el\'ements, et $\Gamma$ n'est donc pas un groupe cci. 
\hfill $\square$
\enddemo

\bigskip
\head
Vieux pense-b\^ete du \S~XII
\endhead
\medskip

Au cas o\`u il y ait une suite, d\'efinir $C^*_{\pi}(\Gamma)$ !!!

\proclaim{A voir !!!} Pour un groupe $\Gamma$, les deux propri\'et\'es suivantes
sont \'equivalentes~:
\roster
\item"(i)" $\Gamma$ est cci~;
\item"(ii)" pour qu'une repr\'esentation unitaire irr\'eductible
$\pi$ de $\Gamma$ soit de classe
$\Cal C_0$, il faut et il suffit que $\lambda _{\Gamma} \prec \pi$.
\endroster
\endproclaim

Message du Bachir du 9/11/04

Il y a un contre-exemple "stupide":
Prenons $\Gamma$
possedant un sous-groupe distingue $N$
fini et non trivial. Alors la reguliere de
$\Gamma/N$, vue comme representation
de $\Gamma$, est $\Cal C_0$ mais ne contient pas faiblement
la reguliere de $\Gamma.$
On peut egalement  (dans quelques
cas) faire la meme chose avec une
rep irreductible de $\Gamma/N$.
Il faut donc modifier en consequence ta question
et supposer de plus (ce qui est de toute facon
necessaire) que $\pi$ est fidele.
Je ne connais pas la reponse, mais la question
me semble interessante.
Amities,
Bachir

\medskip

\proclaim{Reformulation ???????????????} Pour un groupe $\Gamma$, les deux
propri\'et\'es suivantes sont \'equivalentes~:
\roster
\item"(i)" $\Gamma$ est cci~;
\item"(ii)" pour qu'une repr\'esentation unitaire irr\'eductible
$\pi$ de $\Gamma$ soit de classe
$\Cal C_0$, il faut et il suffit que la C$^*$-alg\`ebre r\'eguli\`ere
$C^*_{\lambda}(\Gamma)$ soit un quotient de la C$^*$-alg\`ebre
$C^*_{\pi}(\Gamma)$.
\endroster
\endproclaim

\bigskip
\head{\bf
$\infty \hdots \infty$.~C$^*$-alg\`ebres r\'eduites des $3$-groupes infinis
}
\endhead
\medskip

Quelques ingr\'edients d'un futur paragraphe ???. 
\newline
??????? Y repenser !!!!!!!
\newline

\medskip

La C$^*$-alg\`ebre r\'eduite d'un groupe poss\'edant un sous-groupe normal
moyennable non r\'eduit \`a $1$ n'est pas simple.
(Pour ce fait bien connu et pour d'autres g\'en\'eralit\'es concernant
la simplicit\'e des C$^*$-alg\`ebres de groupes, voir par exemple \cite{Harp}.)
En particulier, soit $\Gamma$ 
un $3$-groupe infini 
qui est le groupe fondamental d'une vari\'et\'e de Seifert
ou qui est un sous-groupe du groupe $Sol$ 
(le groupe de Lie r\'esoluble de dimension $3$, connexe et simplement connexe)~;
alors la C$^*$-alg\`ebre $C^*_{\lambda}(\Gamma)$ n'est pas simple.

\smallskip

Soient $M$ une vari\'et\'e r\'eductible dont le groupe est un produit libre 
de groupes non r\'eduits \`a un \'el\'ement,
$\Gamma = \Gamma_1 \ast \cdots \ast \Gamma_p$,
avec $p \ge 3$, ou avec
$p = 2$, $\vert \Gamma_1 \vert \ge 3$ et $\vert \Gamma_1 \vert \ge 2$.
Alors la C$^*$-alg\`ebre $C^*_{\lambda}(\Gamma)$ est simple
\cite{PaSa--79}.

\enddocument

\smallskip
\noindent {\it ???i\`eme cas : 
$K \approx \Bbb Z \oplus \Bbb Z$ et $Q$ est infini.}

Si l'indice de $Z(\Delta)$ dans $\Gamma$ est fini,
alors $\Gamma$ est isomorphe \`a un groupe de surface
(th\'eor\`eme 10.6 de \cite{Hemp--76}),
donc au groupe fondamental d'un $2$-tore ou d'une bouteille de Klein
(car $\Gamma$ est cci)~;
alors $\Cal P (M)$ est un $I$-fibr\'e sur etcetcetcetcetcetc
??????? et $\Cal P(M) = M$ ???????

Si l'indice de $Z(\Gamma)$ dans $\Gamma$ est infini,
alors ???????  $\Cal P(M) = M$ ???????
..... (th\'eor\`eme 11.6 de \cite{Hemp--76}).

\smallskip
\noindent {\it Troisi\`eme cas : 
$Z(\Delta) \approx \Bbb Z \oplus \Bbb Z \oplus \Bbb Z$.}
\par

.............

\bigskip

\noindent {\it Bribes d'une ancienne r\'edaction !!!!!}

Supposons donc d\'esormais que $Z(\Delta)$ n'est pas cyclique,
et distinguons deux cas selon que $M$ est orientable ou non.

\smallskip

{\it Cas o\`u $M$ est orientable.} 
Notons $N$ le rev\^etement fini de $M$ 
de groupe fondamental $\Delta$.
Par le th\'eor\`eme de la sph\`ere (chapitre 4 de \cite{Hemp--76}),
la seule obstruction \`a ce que
$N$ soit $\Bbb P^2$-irr\'eductible 
serait que $N$ contienne une \lq\lq boule exotique\rq\rq \ 
(une boule d'homotopie  non hom\'eomorphe \`a une $3$-boule). 

\`A notre connaissance, on ignore s'il est en g\'en\'eral vrai
qu'un rev\^etement fini d'une $3$-vari\'et\'e
ne contenant aucune boule exotique 
ne contient lui non plus aucune boule exotique. 
En revanche, c'est connu dans deux cas particuliers~:
(i) si le rev\' etement est \`a deux feuillets
(lemme 10.4 de \cite{Hemp--76}) et
(ii) si $M$ est une vari\'et\'e de Haken
(corollaire 13.5 \cite{Hemp-76}).

La vari\'et\'e $\Cal P(N)$  est donc $\Bbb P^2$-irr\'eductible. 
Puisque son groupe fondamental $\Delta$ 
est sans torsion et \`a centre non r\'eduit \`a $1$,
il contient un sous-groupe normal cyclique infini.
Il r\'esulte donc de la \lq\lq conjecture\rq\rq \ de Seifert
que $N$ est une vari\'et\'e de Seifert. 
En particulier, le centre $Z(\Delta)$ de son groupe fondamental
est de type fini.

!!!!!!! Suite \`a remettre au net !!!!!!!

Il existe une classification des vari\'et\'es $\Cal P(N)$
possibles au th\'eor\`eme 12.10 de \cite{Hemp--76}~; 
plus pr\'ecisemment $\Cal P(N)$ est 
soit $\Bbb S^1 \times \Bbb S^1 \times I$, 
soit $\Bbb S^1 \times \Bbb S^1 \times \Bbb S^1$. 
Dans le premier cas, $\Gamma$ contient $\Bbb Z\oplus \Bbb Z$ 
comme sous-groupe d'indice fini, 
et avec \cite{Jaco**} $\Cal P(M)$ est lui-m\^eme 
$\Bbb S^1 \times \Bbb S^1 \times I$ 
ou le $I$-fibr\'e non trivial sur la bouteille de Klein. 
En particulier $M$ est un fibr\'e de Seifert, et est Haken.

Nous allons montrer que
$\tilde{M}$ (ou $M$, si elle est orientable) est un fibr\'e de
Seifert, mais aussi que c'est une vari\'et\'e Haken, ce qui
permettra de montrer, lorsque $M$ est non-orientable, que $M$ est
aussi un fibr\'e de Seifert (une 3-vari\'et\'e finiment rev\^ etue
par un fibr\'e de Seifert, Haken, est un fibr\'e de Seifert, cf.
th\'eor\`eme II.6.3, \cite{JaSh--79}).

et avec la conjecture de Seifert, $\Cal
P(\tilde{M})$ est un fibr\'e de Seifert. En particulier $Z(N)$ est
n\'ecessairement de type fini.

Dans le second cas on a la suite exacte :
$$ 
1
\longrightarrow \Bbb Z\oplus\Bbb Z\oplus\Bbb Z
\longrightarrow \Gamma
\longrightarrow G
\longrightarrow 1
$$ 
avec $G$ fini d'ordre $n$. Notons
$A=\Bbb Z\oplus\Bbb Z\oplus\Bbb Z$, et consid\'erons un \.

\bigskip\noindent ??????? \bigskip\noindent

, et il est possible de v\'erifier que chacune d'entre elles est
un fibr\'e de Seifert. En particulier $\pi_1(M) = \pi_1(\Cal P
(M))$ poss\`ede un sous-groupe normal cyclique infini. En
utilisant \`a nouveau la conjecture de Seifert, on s'assure que
$M$ (et non plus seulement $\Cal P(M)$) est un fibr\'e de Seifert.

\smallskip

{\it Cas o\`u $M$ n'est pas orientable.}
Le rev\^etement orientable $\tilde M$  de $M$
est encore $\Bbb P^2$-irr\'eductible (lemme 10.4 de \cite{Hemp--76})
et son groupoe fondamental $\tilde \Gamma$ est encorre cci (lemme 5).

????????????
\hfill $\square$

il est bien connu que $\phi$ est conjugu\'e
\`a une puissance de l'une des matrices
$$
s \, = \,
\left( \matrix 0 & -1 \\ 1 & \phantom{-}1 \endmatrix \right) ,
\hskip.2cm \text{d'ordre $6$,}
\hskip1cm
t \, = \,
\left( \matrix \phantom{-}0 & 1 \\ -1 & 0 \endmatrix \right) ,
\hskip.2cm \text{d'ordre $4$.}
$$
(Cela r\'esulte par exemple de la d\'ecomposition de $SL(2,\Bbb Z)$
en produit libre du groupe cyclique d'ordre $6$ engendr\'e par $s$
et du groupe cyclique d'ordre $4$ engendr\'e par $t$
avec amalgamation sur le groupe d'ordre $2$ engendr\'e par $s^3=t^2$.)
\par\noindent
Alors ??????????? \`a compl\'eter !!!!!!!!!!!!

\demo\nofrills{R\'efl\'echir \`a~: (??????? !!!!!!!)\usualspace}

(iii)  Soit $\Gamma = \Gamma_1 \ast_{\Gamma_0} \Gamma_2$
un produit libre avec amalgamation relativement \`a des inclusions propres
$\Gamma_0 \subset \Gamma_1$ et $\Gamma_0 \subset \Gamma_2$
telles que le seul sous-groupe de $\Gamma_0$ qui est normal
\`a la fois dans $\Gamma_1$ et dans $\Gamma_2$ est le groupe r\'eduit \`a $1$.
Alors $\Gamma$ est cci.
\enddemo

\bigskip

\demo{Ancienne d\'emonstration de la prop. 15} 
La somme des genres des composantes connexes 
du bord d'une $3$-vari\'et\'e est major\'ee par son premier nombre de Betti
(\cite{SeTh--34}, \S \ 64, th\'eor\`eme IV). 
Le premier nombre de Betti de $M$ \'etant nul, 
toutes les composantes connexes du bord $\partial M$ sont des sph\`eres~;
par suite $\hat M$ et $\Cal P (M)$ sont sans bord.
Par ailleurs,  $\Cal P (M)$ est irr\'eductible,
car  $\Bbb Z / 2\Bbb Z$ ne se d\'ecompose pas en
un produit libre non trivial.
Posons d\'esormais $P = \Cal P (M)$~; 
notons $p : \tilde P \longrightarrow P$ son rev\^etement d'orientation,
qui est \`a deux feuillets,
et $\sigma : \tilde P \longrightarrow \tilde P$ 
l'involution correspondante,
qui est sans point fixe et qui pr\'eserve l'orientation.
Remarquons que la $3$-vari\'et\'e $\tilde P$ est ferm\'ee et simplement connexe.

Nous allons montrer que toute $2$-sph\`ere dans $\tilde P$
borde deux $3$-boules, de sorte que $\tilde P$
est hom\'eomorphe \`a $\Bbb S^3$ par un th\'eor\`eme 
\footnote{
Est-ce \lq\lq un th\'eor\`eme d'Alexander \'etablit que
toute sph\`ere plong\'ee dans $S^3$ borde deux boules\rq\rq ,
o\`u un th\'eor\`eme disant qu'une vari\'et\'e
dans laquelle toute sph\`ere de codimension $1$ borde deux boules est une
sph\`ere ?????
}
d'Alexander. 
Comme il n'existe sur $\Bbb S^3$ 
qu'une classe d'involution sans point fixe pr\'eservant l'orientation
\cite{Live--60}, il en r\'esultera que le quotient $P = \Cal P /\sigma$
est bien $\Bbb P^3$.

 Soit $S$ est une $2$-sph\`ere plong\'ee dans $\tilde{P}$.
C'est une surface bilat\`ere, 
car $S$ et $\tilde{M}$ sont orientables
(\cite{SeTh--34}, \S \ 76, th\'eor\`eme III)~;
de plus la sph\`ere $S$ est s\'eparante dans $\tilde P$
par simple connexit\'e de $\tilde P$ 
(autrement, un argument classique permet de montrer l'existence
d'un \'el\'ement d'ordre infini dans $H_1(\tilde{P},\Bbb Z)$
????? invoquer la dualit\'e ad hoc !!!!!!!).
Supposons alors que l'une des composantes connexes de $\tilde P \setminus S$
bord\'ee par $S$ ne soit pas une boule~; 
nous allons achever la preuve en aboutissant \`a une contradiction.

\medskip

\noindent Esquisse de la suite

(1) Etape p\'enible. On se ram\`eme \`a l'une des deux situations suivantes~:
\par\noindent
ou bien $\sigma(S) = S$, ou bien $\sigma(S) \cap S = \emptyset$. Voir~:

(2) Monter alors que $p(S)$ est une $2$-sph\`ere dans $P$
qui est essentielle (???). !!! NON -- $\Bbb P^2$

\medskip\noindent
Suite = ancienne version !!! (sauf que j'ai chang\'e $M$ en $P$).
\medskip

Supposons tout d'abord que $\sigma$ ne pr\'eserve pas $S$, i.e.
que $\sigma(S)\not= S$. Si $\sigma(S)\cap S\not= \emptyset$, alors
en d\'eformant l\'eg\`erement $S$ par isotopie, on peut mettre $S$
en position g\'en\'erale et supposer que $\sigma(S)\cap S$ est une
$1$-vari\'et\'e ferm\'ee $C$, c'est \`a dire une r\'eunion
disjointe d'un nombre fini de cercles.

Puisque $\sigma$ est involutive, $\sigma$ doit pr\'eserver
(globalement) $C$, et de plus $\sigma$ envoie
hom\'eo\-morphiquement les composantes de $S-C$ sur les
composantes de $\sigma(S)-C$. Par finitude du nombre de
composantes, on peut trouver un cercle $C_1\subset C$ qui borde
une r\'egion $D$ dans $S$, qui est un disque ouvert dans $S-C$.
Par construction $\sigma(D)$ est aussi un disque ouvert de
$\sigma(S)-C$, bord\'e par $\sigma(C_1)$.

Si $\sigma(C_1)=C_1$, alors $D\cup C_1\cup\sigma(D)$ est une
sph\`ere plong\'ee pr\'eserv\'ee par $\sigma$. Elle borde deux
r\'egions disjointes $D_1$ et $D_1'$ dans $\tilde{P}$, et l'une
pr\'ecisemment de ces deux r\'egions -- disons $D_1$ -- a son
int\'erieur dans $\tilde{P}-S$. Si $D_1$ est une boule, alors on
peut d\'eformer $S$ pour supprimer l'intersection $C_1$, et
recommencer~; le nombre de composantes de $C$ diminue alors
strictement.

Si $\sigma(C_1)\not=C_1$, $C_1$ borde deux disques dans
$\sigma(S)$, dont un exactement contient $\sigma(C_1)$~; on note
$D'$ celui de ces disques qui ne contient pas $\sigma(C_1)$. Alors
$S_1=D\cup C_1\cup D'$ et $S_2=\sigma(S_1)$ sont deux sph\`eres
plong\'ees. Seulement elles ne sont pas \`a priori disjointes~;
ceci car $D'$ peut intersecter en son int\'erieur la surface $S$,
disons en $C_2,C_3,\ldots ,C_i$. Dans ce cas, si aucune
intersection n'a lieu avec $\sigma(D')$ (c'est \`a dire lorsque
$\sigma(D')$ ne contient aucun des $C_2,C_3,\ldots ,C_i$), les
deux sph\`eres $S_1$ et $S_2$ ont une intersection vide, et on
peut passer \`a l'argument suivant. Sinon on consid\`ere une
composante $C_i$ de $C$ dans $D'$, qui borde un disque ouvert $U$
dans $S-C$, et on applique le m\^eme argument que ci-dessus pour
construire deux autres sph\`eres plong\'ees (car
$\sigma(C_i)\not=C_i$) $S_3,S_4$, mais comme pr\'ec\'edemment
elles ne sont pas forc\'ement disjointes : elles peuvent chacune
contenir un m\^eme nombre de composantes de $C$. Tant que l'on n'a
pas trouv\'e deux sph\`eres disjointes on poursuit ce proc\'ed\'e.
Il finira n\'ecessairement par s'arr\^eter : ceci \`a cause de la
fa\c con dont nous avons choisi $D'$, et du fait que les
composantes de $C$ ne peuvent appara\^\i tre qu'au plus une fois
sur les sph\`eres successivement cr\'ees, en effet, les surfaces
\'etant en position g\'en\'erale, on n'a que des croisements
doubles, et $C_i$ (par exemple) n'appara\^\i t qu'une fois sur $S$
et une fois sur $\sigma(S)$. Par finitude du nombre de composantes
de $C$, on finit par trouver deux sph\`eres, disons $S_1,S_2$,
plong\'ees et disjointes, telles que $\sigma(S_1)=S_2$.

Puisque $S_1,S_2$ sont toutes deux s\'eparantes,
$\tilde{P}-(S_1\cup S_2)$ a trois composantes connexes qui sont :
$E$ dont le bord est $S_1\cup S_2$, $I_1$ dont le bord est $S_1$,
et $I_2$ dont le bord est $S_2$. Alors n\'ecessairement $\sigma$
envoie $I_1$ hom\'eomorphiquement sur $I_2$ et r\'eciproquement.
Si $I_1$ et $I_2$ sont des boules, on peut d\'eformer $S$ par
isotopie pour supprimer les intersections $C_1$ et $\sigma(C_1)$.
Comme pr\'ec\'edemment, le nombre de composantes de $C$ diminue
strictement.

En poursuivant ce proc\'ed\'e, on finit par obtenir, soit :

\noindent
 -- Une sph\`ere $S$, pr\'eserv\'ee par $\sigma$, qui ne
borde pas deux boules.

\noindent
 -- Deux sph\`eres disjointes $S_1$ et
$S_2$, images l'une de l'autre par $\sigma$, chacune s\'eparant
l'espace en deux composantes dont aucune n'est une boule (deux
'fausses' boules).

Dans le dernier cas, $p(S_1)=p(S_2)$ est un sph\`ere essentielle
dans $\Cal P(P)$, ce qui contredit l'irr\'educibilit\'e de $\Cal
P(P)$.

Dans le premier cas $p(S)$ est un plan projectif. Sur $S$,
$\sigma$ renverse l'orientation, et par orientabilit\'e doit donc
permuter les deux composantes connexes de $\tilde{P}-S$, que nous
appelerons $D_1$ et $D_2$. En particulier $D_1$ et $D_2$ sont
hom\'eomorphes, et aucune n'est une boule. Choisissons un
voisinage r\'egulier de $S$ dans $\tilde{P}$, qui soit
pr\'eserv\'e par $\sigma$. C'est un fibr\'e trivial en
intervalles, et nous noterons $S_1$ et $S_2$ ses deux composantes
au bord. Par construction $p(S_1)=p(S_2)$ est une sph\`ere
essentielle dans $\Cal P(P)$, ce qui contredit
l'irr\'educibilit\'e de $\Cal P(P)$. \hfill $\square$
\enddemo

\bigskip

\demo{Ancienne d\'emonstration de la proposition 18} 
Un tel groupe fondamental est un r\'eseau
dans la composante connexe du groupe des isom\'etries de l'espace hyperbolique.
Or c'est une cons\'equence du
th\'eor\`eme de densit\'e de Borel que, 
plus g\'en\'eralement, 
dans un groupe de Lie $G$ connexe simple non compact 
\`a centre r\'eduit \`a un \'el\'ement,
tout r\'eseau $\Gamma$ est cci.

Plus pr\'ecis\'ement, soit $\gamma \in \Gamma$ un \'el\'ement \`a 
classe de conjugaison finie,
c'est-\`a-dire un \'el\'ement tel que
le centralisateur $Z_{\Gamma}(\gamma)$ de $\gamma$ dans $\Gamma$
est d'indice fini dans $\Gamma$~;
alors l'adh\'erence de Zariski $\overline{Z_{\Gamma}(\gamma)}^Z$
est d'indice fini dans l'adh\'erence de Zariski $\overline{\Gamma}^Z$~;
mais $\overline{\Gamma}^Z = G$ par densit\'e de Borel,
donc $\overline{Z_{\Gamma}(\gamma)}^Z = G$ par Zariski-connexit\'e de $G$,
de sorte que $\gamma$ est central dans $G$, et finalement $\gamma = 1$.
\hfill $\square$
\enddemo

\bigskip

Vieux passage de la preuve du lemme 10.

{\it Si $\phi$ est parabolique,} alors $\phi$ est conjugu\'e 
dans $SL(2,\Bbb Z)$ \`a une matrice de la forme
$\pm \left( \matrix 1 & k \\ 0 & 1 \endmatrix \right)$, 
avec $\vert k \vert \ge 1$.
L'automorphisme $\phi$ pr\'eserve donc
un sous-groupe cyclique infini~$C$ de $\Bbb Z^2$,
et ce sous-groupe est normal dans $\pi_1(N)$.
La vari\'et\'e $N$ \'etant de Haken
\footnote{
Rappelons qu'une $3$-vari\'et\'e est dite {\it de Haken}
si elle est irr\'eductible 
et si elle contient une surface bilat\`ere incompressible.
}
\hskip-.1cm
, 
$N$~est une vari\'et\'e de Seifert par 
la \lq\lq conjecture\rq\rq \ des fibr\'es de Seifert
(qui est un r\'esultat nettement plus ancien
dans le cas Haken que dans le cas g\'en\'eral~;
voir le corollaire II.6.4 de \cite{JaSh--79}).

\bigskip
\enddocument